\documentclass[10pt,a4paper,oneside,reqno]{amsart}
\usepackage{amsfonts}
\usepackage{amsthm}
\usepackage{amsmath}
\usepackage{amscd}
\usepackage[latin2]{inputenc}
\usepackage{t1enc}
\usepackage[mathscr]{eucal}
\usepackage{indentfirst}
\usepackage{graphicx}
\usepackage{graphics}
\usepackage{pict2e}
\usepackage{hyperref}
\usepackage{epic}
\numberwithin{equation}{section}
\usepackage[margin=2.9cm]{geometry}
\usepackage{epstopdf}
\usepackage{enumitem}

\theoremstyle{plain}

 \theoremstyle{definition}

\newtheorem{?}[Th]{Problem}

\begin{document}

\title{ High-Frequency Instabilities of the Kawahara Equation: 
 A Perturbative Approach}

\author{Ryan Creedon$^{1}$, Bernard Deconinck$^2$, Olga Trichtchenko$^3$ \\~\\ \fontsize{0.15}{0.15}\selectfont $^1$ Dept. of Applied Mathematics, U. of Washington,  Seattle, WA, 98105, USA (\emph{creedon@uw.edu}) \\ $^2$ Dept. of Applied Mathematics, U. of Washington,  Seattle, WA, 98105, USA (\emph{deconinc@uw.edu}) \\$^3$ Dept. of Physics and Astronomy, U. of Western Ontario, London, ON, N6A 3K7, CA (\emph{otrichtc@uwo.ca}) \footnotesize \\ ~ \\ January 17, 2021}





\begin{abstract} We analyze the spectral stability of small-amplitude, periodic, traveling-wave solutions of the Kawahara equation. These solutions exhibit high-frequency instabilities when subject to bounded perturbations on the whole real line. We introduce a formal perturbation method to determine the asymptotic growth rates of these instabilities, among other properties. Explicit numerical computations are used to verify our asymptotic results.
\end{abstract}

\maketitle

\section{Introduction} 
\begin{figure}[t!]
\centering \includegraphics[width=15cm,height=7cm]{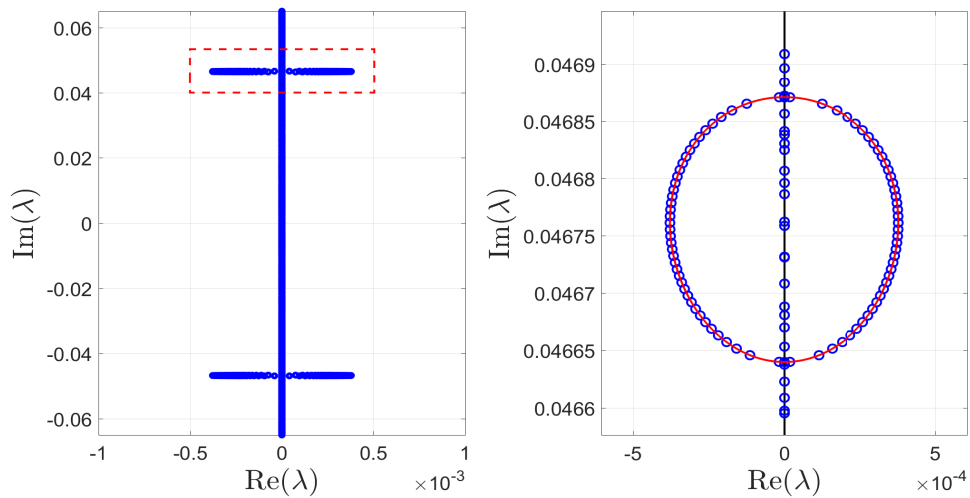}
\caption{(Left) A stability spectrum of Stokes wave solutions of \eqref{1} with $\alpha = 1$, $\beta = 0.7$, $\sigma = 1$, and small-amplitude parameter $\varepsilon = 10^{-3}$, computed using the FFH method. A uniform grid of $10^3$ Floquet exponents between $[-1/2,1/2]$ approximates purely imaginary point spectra but misses the high-frequency isolas. A uniform grid of $4\times 10^3$ Floquet exponents in the interval described by \eqref{27b}, obtained in Section 4, captures these isolas. (Right) Zoom-in of the high-frequency isola boxed in the left plot (with fewer point spectra shown for ease of visibility). The red curve is obtained in Section 5 and approximates the isola.   }
\label{fig1}
\end{figure}

We investigate small-amplitude, $L$-periodic, traveling-wave solutions of the Kawahara equation
\begin{align} \label{1}
u_t = \alpha u_{xxx} + \beta u_{5x} + \sigma (u^2)_x,
\end{align}
where $\alpha$, $\beta$, and $\sigma$ are nonzero, real parameters \cite{kawahara72}. Similar to Stokes waves of the Euler water wave problem \cite{stokes1847,whitham67}, these solutions are obtained order by order as power series in a small parameter that scales with the amplitude of the solutions; see \cite{haragus06} and Section 2 below for more details. We refer to these solutions as the Stokes waves of the Kawahara equation.   \\\\
The Kawahara equation is dispersive with linear dispersion relation
\begin{align}
\omega(k) = \alpha k^3 - \beta k^5.
\end{align}
The equation is Hamiltonian,
\begin{align}
u_t = \partial_x \frac{\delta \mathcal{H}}{\delta u},
\end{align}
with
\begin{align} \label{H}
\mathcal{H} = \int_0^L \left(  -\frac{\alpha}{2} u_x^2 + \frac{\beta}{2} u_{xx}^2 + \frac{\sigma}{3} u^3 \right) d x.
\end{align}
In an appropriate traveling frame, the Stokes waves of \eqref{1} are critical points of the Hamiltonian, prompting an investigation of the flow generated by \eqref{H} about the Stokes wave solutions.
 \\\\
Perturbing the Stokes waves by functions bounded in space and exponential in time yields a spectral problem whose spectral elements characterize the temporal growth rates of the perturbations; see Section 3 for more details. We refer to this collection of spectral elements as the stability spectrum of the Stokes waves. \\\\ A standard argument \cite{kapitulapromislow13} shows that the stability spectrum is purely continuous, but Floquet theory can decompose the spectrum into an uncountably infinite collection of point spectra. Each point spectra is indexed by a real number, called the Floquet exponent, that is contained within a compact interval of the real line \cite{haraguskapitula08,johnson10}. \\\\
For the Euler water wave problem, these point spectra depend analytically on the amplitude of the Stokes waves \cite{akersnicholls12,akersnicholls14}. Based on numerical experiments \cite{trichtchenko18}, similar results appear to hold for the Kawahara equation. The spectrum also exhibits quadrafold symmetry due to the underlying Hamiltonian nature of \eqref{1} \cite{haraguskapitula08},\cite{mackay86}. Therefore, for a Stokes wave with given amplitude to be spectrally stable, all point spectra must be on the imaginary axis. Otherwise, there exist perturbations to the Stokes waves that grow exponentially in time, and the Stokes waves are spectrally unstable. \\\\
In contrast with the completely integrable KdV equation ($\beta = 0$) \cite{deconinckkapitula10,nivala10,pava14}, considerably less is known about the stability spectrum of Stokes waves to \eqref{1}. Haragus, Lombardi, and Scheel \cite{haragus06} prove that this spectrum lies on the imaginary axis for small-amplitude Stokes waves in a particular scaling regime. Such solutions are, therefore, spectrally stable. Work by Trichtchenko, Deconinck, and Koll\'{a}r \cite{trichtchenko18} develops necessary criteria for the stability spectrum of a broader class of small-amplitude Stokes waves to leave the imaginary axis and provide numerical evidence for the high-frequency instabilities that result. \\\\
High-frequency instabilities arise from pairwise collisions of nonzero, imaginary elements of the stability spectrum. Upon colliding, these elements may symmetrically bifurcate from the imaginary axis as the amplitude of the Stokes wave grows, resulting in instability \cite{deconinck17},\cite{mackay86}. An example of a high-frequency instability for a small-amplitude Stokes wave of \eqref{1} is seen in \autoref{fig1}. We refer to the locus of spectral elements off the imaginary axis and bounded away from the origin as high-frequency isolas. The isolas of \autoref{fig1}, as well as the rest of the stability spectrum, are obtained numerically using the Floquet-Fourier-Hill (FFH) method; see \cite{deconinck06} for a detailed description of this method.
 ~\\\\
\noindent High-frequency instabilities are not as well-studied as the modulational (or Benjamin-Feir) instability that arises from collisions of spectral elements at the origin of the complex spectral plane \cite{benjamin67},\cite{bridgesmielke95}. Current understanding of high-frequency instabilities is limited mostly to numerical experiments. Exceptions include the works of Akers \cite{akers15} and Trichtchenko, Deconinck, and Koll\'{a}r \cite{trichtchenko18}, which obtain asymptotic information about the high-frequency isolas for the Euler problem in infinitely deep water and for the Kawahara equation, respectively.  \\\\
The purpose of our present work is to build on these results. In particular, for sufficiently small-amplitude solutions, we seek the following:
\begin{enumerate}[label=(\roman*)]
\item the asymptotic range of Floquet exponents that parameterize the high-frequency isolas observed in numerical computations of the stability spectrum,
\item  asymptotic estimates of the most unstable spectral elements of the high-frequency isolas, and
\item expressions for curves asymptotic to these isolas, as seen in \autoref{fig1}.
\end{enumerate}
To obtain these quantities, we develop a perturbation method inspired by \cite{akers15} that readily extends to higher-order calculations. Asymptotic results obtained by this method are then compared with numerical results from the FFH method.


\section{Small-Amplitude Stokes Waves}
We move to a frame traveling with velocity $c$  so that $x \rightarrow x - ct$. Equation \eqref{1} becomes
\begin{align} \label{2}
u_t = cu_x + \alpha u_{xxx} + \beta u_{5x} + \sigma (u^2)_x.
\end{align}
We seek $L$-periodic, steady-state solutions of \eqref{2}. Equating time derivatives to zero and integrating in $x$, we arrive at
\begin{align} \label{3}
cu + \alpha u_{xx} + \beta u_{4x} + \sigma u^2 = \mathcal{C},
\end{align}
where $\mathcal{C}$ is a constant of integration. Using the Galilean symmetry of \eqref{1}, there exists a boost $\xi$ such that, with $c \rightarrow c + \xi$ and $u \rightarrow u + \xi$, $\mathcal{C}$ can be omitted from \eqref{3}:
\begin{align} \label{4}
cu + \alpha u_{xx} + \beta u_{4x} + \sigma u^2 = 0.
\end{align} Scaling $x \rightarrow 2\pi x/L$ and $u \rightarrow 2\pi u/(\alpha L)$ allows us to consider $2\pi$-periodic solutions of 
\begin{align} \label{5}
cu + u_{xx} + \beta u_{4x} + \sigma u^2 = 0,
\end{align}
without loss of generality, provided $c$, $\beta$, and $\sigma$ are appropriately redefined.\\\\
Let $u=u_S(x;\varepsilon)$ be a one-parameter family of $2\pi$-periodic solutions of \eqref{5} with corresponding velocity $c=c(\varepsilon)$. The existence of such a family is rigorously justified by Lyapunov-Schmidt reduction; see \cite{haragus06}. In what follows, we define the parameter $\varepsilon$ as twice the first Fourier coefficient of $u_S(x;\varepsilon)$:
\begin{align} \label{eps_def}
\varepsilon := 2 \widehat{u_{S}(x;\varepsilon)}_1 = \frac{1}{\pi} \int_{0}^{2\pi} u_S(x;\varepsilon) e^{ix} d x,
\end{align}
where $~\widehat{\cdot}~$ is the Fourier transform on the interval $(0,2\pi)$. Because the $\textrm{L}^2(0,2\pi)$ norm of $u_S(x;\varepsilon)$ scales like  $\varepsilon$ when $|\varepsilon| \ll 1$, we call $\varepsilon$ the small-amplitude parameter. \\\\
From \cite{haragus06}, expansions for $u_S(x;\varepsilon)$ and $c(\varepsilon)$ take the form
\begin{subequations} \label{6}
\begin{align}
u_S(x;\varepsilon) &= \sum_{k=1}^{\infty} u_k(x) \varepsilon^k, \\
c(\varepsilon) &= \sum_{k=0}^{\infty} c_{2k} \varepsilon^{2k},
\end{align}
\end{subequations}
where  $u_k(x)$ is analytic and $2\pi$-periodic for each $k$. Exploiting the invariance of \eqref{5} under $x \rightarrow -x$ and $x \rightarrow x + \phi$, we require $u_k(x) = u_k(-x)$ so that $u_S(x;\varepsilon)$ is even in $x$ without loss of generality.
Substituting these expansions into \eqref{5} and following a Poincar\'{e}-Lindstedt perturbation method \cite{whitham67}, one finds corrections to $u_S(x;\varepsilon)$ and $c(\varepsilon)$ order by order. \\\\ One difficulty occurs at leading order of the Poincar\'{e}-Lindstedt method. Substituting expansions \eqref{6} into \eqref{5} and collecting terms of $\mathcal{O}(\varepsilon)$, we find
\begin{align} \label{7}
\left[c_0 + \partial_x^2 + \beta \partial_x^4 \right] u_1(x) = 0.
\end{align}
From \eqref{eps_def}, $\widehat{u_1(x)}_1 = 1/2$. Taking the Fourier transform of \eqref{7} and evaluating at the first mode, we find
\begin{align}
\left[c_0 - 1 + \beta \right]\widehat{u_1(x)}_1 = \frac12(c_0 - 1 + \beta) = 0,
\end{align}
which implies that
\begin{align} \label{8}
c_0 = 1 - \beta.
\end{align}
By inspection,
\begin{align} \label{9}
u_1(x) = \cos(x)
\end{align}
is a solution to \eqref{7} that is analytic, $2\pi$-periodic, even in $x$, and satisfies the normalization $\widehat{u_1(x)}_1 = 1/2$. If $\beta = 1/(1+N^2)$ for any integer $N > 1$, then
\begin{align}
u_1(x) = \cos(x) + C_N \cos(Nx),
\end{align}
where $C_N$ is an arbitrary real constant, is an equally valid solution to \eqref{7} with the requisite properties. In this case, the Stokes waves are said to be resonant and exhibit Wilton ripples \cite{wilton1915}. Expansions \eqref{6} must be modified as a result; see \cite{akers12,haupt88}, for instance. \\\\
In this manuscript, we restrict to nonresonant Stokes waves:
\begin{align} \label{10}
\beta \neq \frac{1}{1+N^2},
\end{align}
for $N$ stated above, and  \eqref{8} and \eqref{9} are the unique leading-order behaviors of $c(\varepsilon)$ and $u_S(x;\varepsilon)$, respectively. The remainder of the Poincar\'{e}-Lindstedt method follows as usual. We terminate the method after third-order corrections, as this is sufficient for our calculations that follow. We find \\
\begin{subequations} \label{11}
\begin{align}
u_S(x;\varepsilon) &=  \varepsilon u_{1}(x) + \varepsilon^2 u_{2}(x) + \varepsilon^3 u_{3}(x) + \mathcal{O}(\varepsilon^4) \\
&= \varepsilon \cos(x) + \varepsilon^2 \frac{\sigma}{2} \left( -\frac{1}{c_0} + \frac{2}{\Omega(2)} \cos(2x) \right) + \varepsilon^3\frac{3\sigma^2}{\Omega(2)\Omega(3)}\cos(3x)  + \mathcal{O}(\varepsilon^4) \nonumber, \\
c(\varepsilon) &= c_0 + c_2 \varepsilon^2 + \mathcal{O}(\varepsilon^4) \\
&= 1-\beta + \sigma^2 \left(\frac{1}{c_0} -\frac{1}{\Omega(2)}\right)\varepsilon^2  +\mathcal{O}(\varepsilon^4), \nonumber
\end{align}
\end{subequations} ~\\
where $\Omega(\cdot)$ is the linear dispersion relation of the Kawahara equation \eqref{1} (with $\alpha = 1$) in a frame traveling at velocity $c(\varepsilon)$:
\begin{align} \label{12}
\Omega(k) = -c_0k + k^3 - \beta k^5.
\end{align}
\section{Necessary Conditions for High-Frequency Instability} \vspace*{0.2cm}
\subsection{ The Stability Spectrum \nopunct} ~\\\\
\indent We consider a perturbation to $u_S(x;\varepsilon)$ of the form
\begin{align}
u(x,t) = u_S(x;\varepsilon) + \rho v(x,t) + \mathcal{O}(\rho^2),
\end{align}
where $|\rho| \ll 1$ is a small parameter independent of $\varepsilon$ and $v(x,t)$ is a sufficiently smooth, bounded function of $x$ on the whole real line for each $t\geq 0$. Substituting \eqref{1} (with $\alpha = 1$) and collecting terms of $\mathcal{O}(\rho)$,  we find by formally separating variables
\begin{align}
v(x,t) = e^{\lambda t} W(x) + c.c.,
\end{align}
where  $c.c.$ denotes complex conjugation of what precedes and $W(x)$ satisfies the spectral problem
\begin{align} \label{13}
\lambda W(x) = \mathcal{L}(u_S(x;\varepsilon),c(\varepsilon),\beta,\sigma) W(x),
\end{align}
with
\begin{align}
\mathcal{L}(u_S(x;\varepsilon),c(\varepsilon),\beta,\sigma)= c(\varepsilon)\partial_x + \partial_x^3 + \beta \partial_x^5 + 2\sigma u_S(x;\varepsilon) \partial_x + 2\sigma u_S'(x;\varepsilon),
\end{align}
where primes denote differentiation with respect to $x$. From Floquet theory \cite{haraguskapitula08}, all solutions of \eqref{13} that are bounded over $\mathbb{R}$ take the form
\begin{align} \label{14}
W(x) = e^{i\mu x} w(x),
\end{align}
where $\mu \in [-1/2,1/2]$ is the Floquet exponent and $w(x)$ is $2\pi$-periodic in an appropriately chosen function space. \\\\
\textbf{Remark}. The conjugate of $W(x)$ is a solution of \eqref{13} with spectral parameter $\overline{\lambda}$. Since the spectrum of $\mathcal{L}$ is invariant under conjugation according to \cite{haraguskapitula08}, one can restrict $\mu$ to the interval $[0, 1/2]$ without loss of generality. \\\\
 Substituting \eqref{14} into \eqref{13}, our spectral problem becomes a one-parameter family of spectral problems:
\begin{align} \label{15}
\lambda^{\mu} w(x) = \mathcal{L}^{\mu}(u_S(x;\varepsilon),c(\varepsilon),\beta,\sigma) w(x),
\end{align}
where $\mathcal{L}^{\mu}$ is $\mathcal{L}$ with $\partial_x \rightarrow i\mu + \partial_x$. In light of \eqref{15}, we require $w(x) \in \textrm{H}^5_{\textrm{per}}(0,2\pi)$ so that $\mathcal{L}^{\mu}$ is a closed operator densely defined on the separable Hilbert space $\textrm{L}_{\textrm{per}}^2(0,2\pi)$ for a given $\mu$. Then, $\mathcal{L}^\mu$ has a discrete spectrum of eigenvalues $\lambda^{\mu}$ for each $\mu$ and the union of $\lambda^{\mu}$ over all $\mu \in [0,1/2]$ yields the purely continuous spectrum of $\mathcal{L}$, which is the stability spectrum of Stokes waves. See \cite{haraguskapitula08} for more details.\\\\
As stated in the introduction, if there exists $\lambda^\mu$ with $\textrm{Re}\left(\lambda^\mu\right)>0$, then there exists a perturbation to the Stokes wave that grows exponentially in time, and we say that the Stokes wave is spectrally unstable. Otherwise, the wave is spectrally stable. Since \eqref{15} is obtained from a linearization of a Hamiltonian system \eqref{1}, the stability spectrum is invariant under conjugation and negation. As a result, spectral stability implies that all eigenvalues of $\mathcal{L}^\mu$ are purely imaginary. \\
\subsection{ The Necessary Conditions for High-Frequency Instability \nopunct} ~\\\\
\indent For fixed $\mu$, the operator $\mathcal{L}^\mu$ depends implicitly on the small parameter $\varepsilon$ through its dependence on $u_S(x;\varepsilon)$ and $c(\varepsilon)$. If $\varepsilon = 0$, $\mathcal{L}^\mu$ reduces to
\begin{align} \label{16.8}
\mathcal{L}^{\mu}_0 = c_0(i\mu+\partial_x) + (i\mu+\partial_x)^3 + \beta (i\mu+\partial_x)^5,
\end{align}
 a constant-coefficient operator, and its eigenvalues $\lambda^{\mu}_0$ are explicitly given by
 \begin{align}  \label{16} \lambda^{\mu}_{0,n} &= -i \Omega(\mu+n),
\end{align}
where $n \in \mathbb{Z}$. For all $\mu$ and $n$, these eigenvalues are purely imaginary, implying that the zero-amplitude solution of the Kawahara equation is spectrally stable. \\\\
Importantly, not all eigenvalues given by \eqref{16} are simple. Using the theory outlined in \cite{kollar19} and \cite{trichtchenko18}, one has \\\\
\textbf{Theorem 1}. For each $\Delta n \in \mathbb{N}$, there exists a unique Floquet exponent $\mu_0 \in [0, 1/2]$ and unique integers $m$ and $n$ such that $m - n = \Delta n$ and
\begin{align} \label{17}
\lambda^{\mu_0}_{0,n} = \lambda^{\mu_0}_{0,m} \neq 0,
\end{align}
provided that the parameter $\beta$ is nonresonant \eqref{10} and that $\beta$ satisfies the inequality\footnote{A similar statement holds for $\Delta n < 0$. This yields the complex conjugate eigenvalues that satisfy \eqref{17}.}:
\begin{subequations} \label{18}
\begin{align}
\textrm{max} \left(\frac{3}{5(\Delta n)^2},\frac{1}{1+(\Delta n)^2} \right) < &~\beta < \textrm{min} \left(\frac{6}{5(\Delta n)^2},\frac{1}{\left(\frac{\Delta n}{2}\right)^2 +1} \right), \quad \Delta n < 3, \\ \frac{1}{1+(\Delta n)^2} < &~\beta < \frac{1}{1+\left( \frac{\Delta n}{2}\right)^2}, \quad \Delta n \geq 3.
\end{align}
\end{subequations} ~\\\\
The proof is found in the Appendix. \\\\
The eigenfunctions of these nonsimple eigenvalues take the form
\begin{align}
w_0(x) = \gamma_{0} e^{imx} +  \gamma_{1} e^{inx},
\end{align}
where $\gamma_{0},\gamma_{1}$ are arbitrary, complex constants. We assume the eigenvalues that satisfy \eqref{17} are semi-simple with geometric and algebraic multiplicity two. Then, these eigenvalues represent the collision of two simple eigenvalues at $\varepsilon = 0$, and \eqref{17} is referred to as the collision condition. \\\\
Collision of eigenvalues away from the origin is a necessary condition for the development of high-frequency instabilities. Inequality \eqref{18} guarantees that there are a finite\footnote{For $\Delta n$ sufficiently large, $\beta$ fails to satisfy inequality \eqref{18}, and no high-frequency instabilities occur.} number of such collisions for a given $\beta$: this is in contrast to the water wave problem, where a countably infinite number of collisions occur \cite{deconinck17}. Each collision site can be enumerated by $\Delta n$. The largest high-frequency isola occurs from the $\Delta n = 1$ collision, which we study in Section 4.1.  \\\\
A second condition for high-frequency instabilities necessitates that the Krein signatures  \cite{krein51} of the two collided eigenvalues have opposite signs \cite{mackay86}. It is shown in \cite{deconinck17, kollar19,trichtchenko18} that this condition is equivalent to
\begin{align} \label{19}
(\mu_0+n)(\mu_0+m) < 0,
\end{align}
where $\mu_0$, $m$, and $n$ are obtained from the collision condition \eqref{17}. For any $\beta$ that satisfies conditions \eqref{10} and \eqref{18} and any $\mu_0$, $m$, and $n$ that satisfies the condition \eqref{17}, \eqref{19} is automatically satisfied; see \cite{kollar19} and \cite{trichtchenko18} for the proof. \\\\
As $|\varepsilon|$ increases in magnitude, a neighborhood of spectral elements around the collided eigenvalues of $\mathcal{L}^{\mu_0}_0$ \eqref{16.8} can leave the imaginary axis, generating high-frequency instabilities. This is seen explicitly in \autoref{fig1} for the parameter choice $\beta = 0.7$, where a $\Delta n = 1$ collision occurs at $\varepsilon = 0$.

\section{Asymptotics of High-Frequency Instabilities} \vspace*{0.2cm}
We obtain spectral data of $\mathcal{L}^\mu$ as a power series expansion in $\varepsilon$ about the collided eigenvalues of $\mathcal{L}^{\mu_0}_0$. First, we apply our method to the largest high-frequency instability corresponding to $\Delta n = 1$. Then, we consider $\Delta n \geq 2$. \\\\
\subsection{High-Frequency Instabilities: $\Delta n = 1$ \nopunct} 
\begin{figure}[t!]
\centering \hspace*{-1cm} \includegraphics[width=15cm,height=7cm]{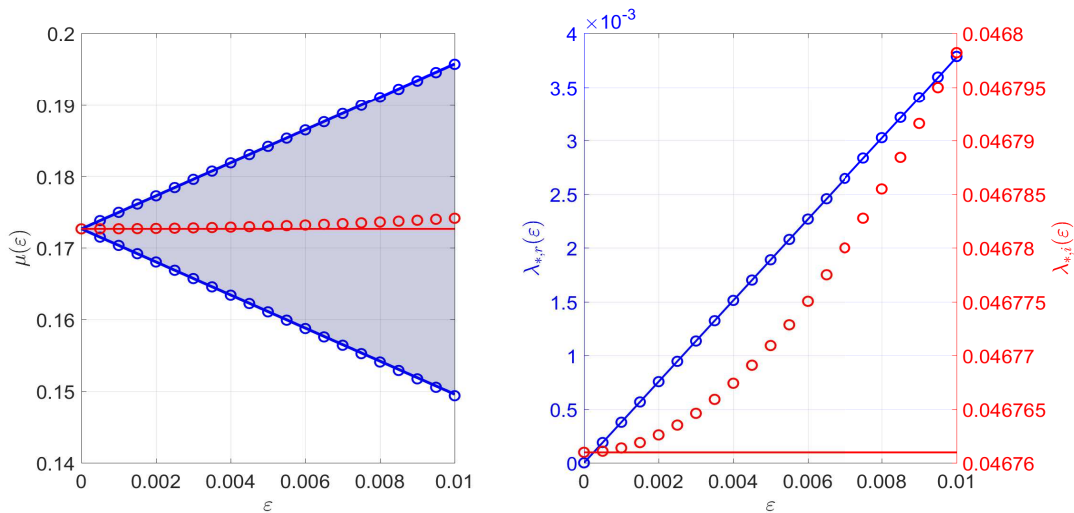}
\caption{(Left) Interval of Floquet exponents that parameterize the $\Delta n = 1$ high-frequency isola for parameters $\alpha = 1$, $\beta = 0.7$, and $\sigma = 1$ as a function of $\varepsilon$. Solid blue curves indicate the asymptotic boundaries of this interval according to \eqref{27b}. Blue circles indicate the numerical boundaries computed using FFH. The solid red curve gives the Floquet exponent corresponding to the most unstable spectral element of the isola according to \eqref{27cc}. The red circles indicate the same but computed numerically using FFH. (Right) The real (blue) and imaginary (red) parts of the most unstable spectral element of the isola as a function of $\varepsilon$. Solid curves illustrate asymptotic result \eqref{27dd}. Circles illustrate results of FFH. }
\label{fig2}
\end{figure}

\begin{figure}[t!]
\centering \hspace*{-1cm} \includegraphics[width=15cm,height=7cm]{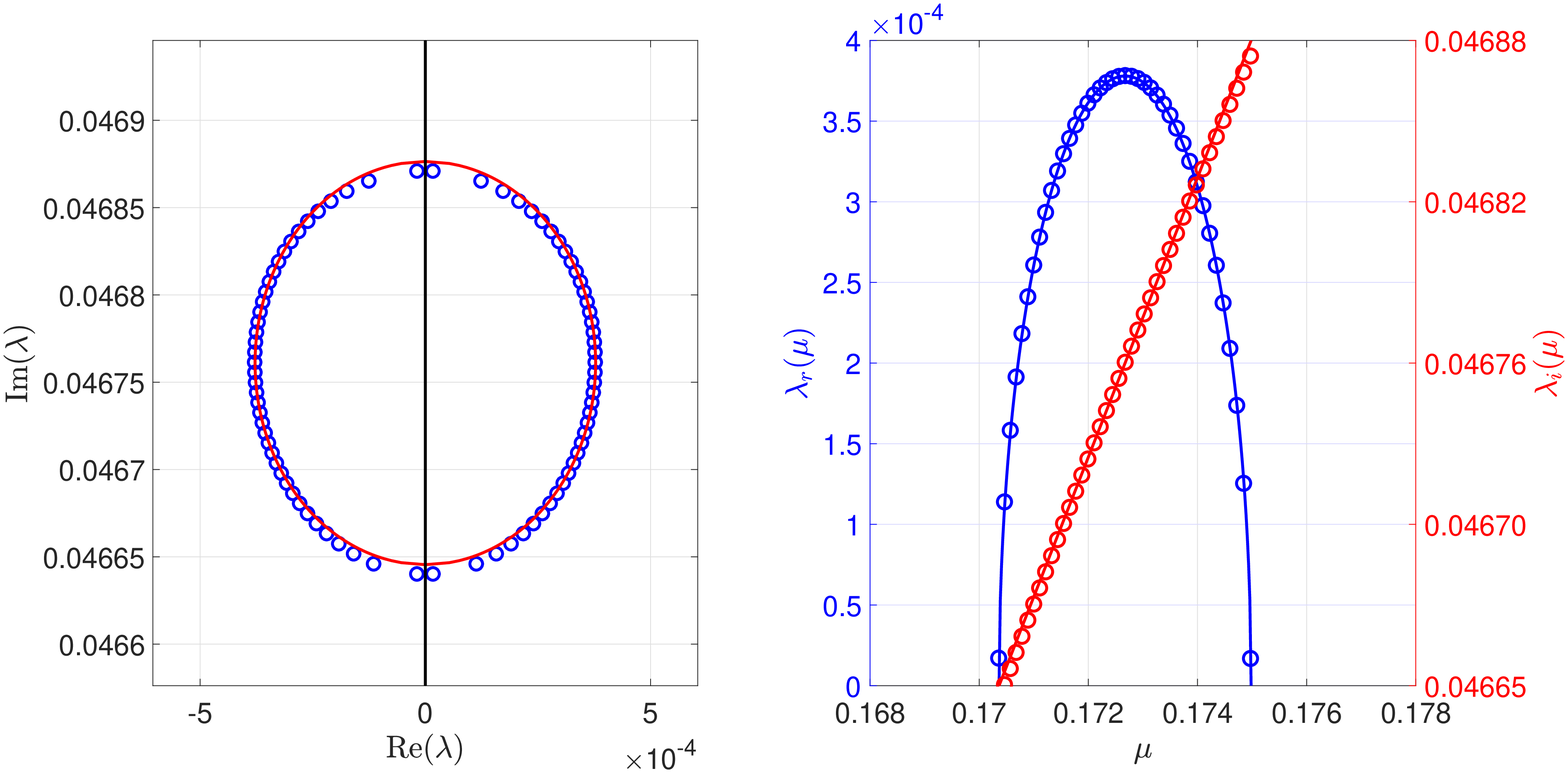}
\caption{(Left) $\Delta n = 1$ high-frequency isola for $\alpha = 1$, $\beta = 0.7$, $\sigma = 1$, and $\varepsilon = 10^{-3}$. The solid red curve is ellipse \eqref{27d}. Blue circles are a subset of spectral elements from the numerically computed isola using FFH. (Right) Floquet parameterization of the real (blue) and imaginary (red) parts of the isola. Solid curves illustrate asymptotic result \eqref{26c}. Circles indicate results of FFH. }
\label{fig3}
\end{figure}

~\\\\
\indent Let $m$ and $n$ be the unique integers that satisfy the collision condition \eqref{17} with $m - n =1$, and let $\mu_0$ be the corresponding unique Floquet exponent in $[0, 1/2]$. Then, the spectral data of $\mathcal{L}^{\mu_0}_0$ that gives rise to a $\Delta n = 1$ high-frequency instability is
\begin{subequations} \label{20}
\begin{align}
\lambda_0 := \lambda^{\mu_0}_{0,n} &=  -i\Omega(\mu_0+n)= \lambda^{\mu_0}_{0,m} = -i\Omega(\mu_0+m) \neq 0,\\
w_0(x) &:= \gamma_0 e^{imx}+ \gamma_1e^{inx}.
\end{align}
\end{subequations}
As $|\varepsilon|$ increases, we assume these data depend analytically on $\varepsilon$:
\begin{subequations}
\begin{align}
\lambda(\varepsilon) &= \lambda_0 + \varepsilon \lambda_1 + \mathcal{O}(\varepsilon^2),  \label{21a} \\
w(x;\varepsilon) &= w_0(x) + \varepsilon w_1(x) + \mathcal{O}(\varepsilon^2), \label{21b}
\end{align}
\end{subequations}
where $\lambda(\varepsilon)$ and $w(x;\varepsilon)$ solve the spectral problem \eqref{15}. \\\\ If $\lambda_0$ is a semi-simple and isolated eigenvalue of $\mathcal{L}^{\mu_0}_0$, \eqref{21a} and \eqref{21b} may be justified using results of analytic perturbation theory, provided $\mu_0$ is fixed \cite{kato66}. Fixing the Floquet exponent in this way, however, gives at most two elements on the high-frequency isola (provided $|\varepsilon|$ is sufficiently small) and these elements do not, in general, correspond to the spectral elements of largest real part on the isola. For these reasons, we expand the Floquet exponent about its resonant value as well:
\begin{align} \label{22}
\mu = \mu(\varepsilon) = \mu_0 + \varepsilon \mu_1 + \mathcal{O}(\varepsilon^2).
\end{align}
As we shall see, $\mu_1$ is constrained to an interval of values that parameterizes an ellipse asymptotic to the high-frequency isola. \\\\
Like Akers \cite{akers15}, we impose the following normalization condition on our eigenfunction $w(x;\varepsilon)$ for uniqueness:
\begin{align} \label{22b}
\widehat{w(x;\varepsilon)}_{n} = 1.
\end{align}
Substituting \eqref{21b} into this normalization condition, we find $\widehat{w_0(x)}_{n} = 1$ and
$\widehat{w_j(x)}_{n} = 0$ for $j \in \mathbb{N}$, meaning $w_0(x)$ fully supports the $n^{\textrm{th}}$ Fourier mode of the eigenfunction $w(x;\varepsilon)$.  As a consequence, \eqref{20} becomes
\begin{align} \label{23}
w_0(x) = e^{inx} + \gamma_0 e^{imx}.
\end{align}
Although $w_0(x)$ does not appear unique at this order, we find an expression for $\gamma_0$ at the next order. \\\\
\textit{The $\mathcal{O}(\varepsilon)$ Problem} ~\\\\
\indent Substituting \eqref{21a}, \eqref{21b}, and \eqref{22} into \eqref{15} and collecting terms of $\mathcal{O}(\varepsilon)$ yields
\begin{align} \label{24}
(\mathcal{L}_0^{\mu_0}-\lambda_0)w_1(x) = \lambda_1w_0(x) - \mathcal{L}_1w_0(x),
\end{align}
where
\begin{align}
\mathcal{L}_1 =  ic_0\mu_1 +3i\mu_1(i\mu_0+\partial_x)^2 + 5i\beta \mu_1 (i\mu_0+\partial_x)^4 + 2\sigma u_1(x)(i\mu_0 + \partial_x) + 2\sigma u_1'(x).
\end{align}
Using \eqref{11} to replace $u_1(x)$, \eqref{23} to replace $w_0(x)$, and $m - n = 1$, \eqref{24} becomes
\begin{align} \label{25}
(\mathcal{L}^{\mu_0}_0-\lambda_0) w_1(x)= &\left[\lambda_1 + i\mu_1 c_g(\mu_0+n)-i\sigma\gamma_{0}(\mu_0+n) \right]e^{inx} \\ &+ \left[\gamma_{0}\left(\lambda_1 +i\mu_1 c_g(\mu_0+m)\right)-i\sigma(\mu_0+m) \right]e^{imx} \nonumber \\ &-i \sigma (\mu_0+n-1)e^{i(n-1)x} -i\sigma \gamma_{0} (\mu_0 +m+1)e^{i(m+1)x}, \nonumber
\end{align}
where $c_g(k) =\Omega'(k)$ is the group velocity of $\Omega$. \\\\
If \eqref{25} can be solved for  $w_1(x) \in \textrm{H}^5_{\textrm{per}}\left(0,2\pi\right)$, the Fredholm alternative necessitates that the inhomogeneous terms on the RHS of \eqref{25} must be orthogonal\footnote{In the $\textrm{L}^2_{\textrm{per}}\left(0,2\pi\right)$ sense} to the nullspace of $(\mathcal{L}^{\mu_0}_0-\lambda_0) ^{\dagger}$, the hermitian adjoint of $\mathcal{L}^{\mu_0}_0-\lambda_0$.  A quick computation shows that $\mathcal{L}^{\mu_0}_0-\lambda_0 $ is skew-Hermitian, and so its nullspace coincides with that of its Hermitian adjoint. The nullspace of $\mathcal{L}^{\mu_0}_0-\lambda_0 $ is, by construction,
\begin{align}
\textrm{Null}(\mathcal{L}^{\mu_0}_0-\lambda_0) = \textrm{Span}\left( e^{inx},e^{imx}\right).
\end{align}
Thus, the solvability conditions for \eqref{25} are
\begin{subequations}
\begin{align}
\left<e^{inx}, \left[\lambda_1 + i\mu_1 c_g(\mu_0+n)-i\sigma\gamma_{0}(\mu_0+n) \right]e^{inx}\right> &= 0  \label{26a},\\
 \left<e^{imx},\left[\gamma_{0}\left(\lambda_1 +i\mu_1 c_g(\mu_0+m)\right)-i\sigma(\mu_0+m) \right]e^{imx} \right> &= 0, \label{26b}
\end{align}
\end{subequations}
where $\left<\cdot,\cdot\right>$ is the standard inner product on $\textrm{L}_{\textrm{per}}^2(0, 2\pi)$. \\\\
\textbf{Remark}. Both solvability conditions can be reinterpreted as removing secular terms from \eqref{25}. Moreover, solvability condition \eqref{26a} coincides with the normalization condition $\widehat{w_1(x)}_n = 0$. \\\\
The solvability conditions \eqref{26a} and \eqref{26b} yield a nonlinear system for $\lambda_1$ and $\gamma_0$ with solution
\begin{subequations}
\begin{align}
\lambda_1 =& -i\mu_1 \left(\frac{c_g(\mu_0+m)+c_g(\mu_0+n)}{2} \right) \label{26c} \\ \quad & \quad \pm \sqrt{-\mu_1^2\left[\frac{c_g(\mu_0+m)-c_g(\mu_0+n)}{2} \right]^2-\sigma^2(\mu_0+m)(\mu_0+n)}, \nonumber\\~ \nonumber\\
\gamma_{0} =& \frac{i\sigma(\mu_0+m)}{\lambda_1 + i\mu_1c_g(\mu_0 +m)}.
\end{align}
\end{subequations}
If $\mu_1 \in (-M_1,M_1)$ with
\begin{align} \label{27a}
M_1 = \frac{2|\sigma| \sqrt{-(\mu_0+m)(\mu_0+n)}}{\left|c_g(\mu_0+m) - c_g(\mu_0+n) \right|},
\end{align}
it follows that $\lambda_1$ has nonzero real part, since $(\mu_0+m)(\mu_0+n
) <0$ by the choice of $\beta$. Therefore, to $\mathcal{O}(\varepsilon)$, the $\Delta n = 1$ high-frequency instability is parameterized by
\begin{align} \label{27b}
\mu \in (\mu_0 -\varepsilon M_1,\mu_0+\varepsilon M_1).
\end{align}
This interval is asymptotically close to the numerically observed interval of Floquet exponents that parameterize the high-frequency isola for $|\varepsilon| \ll 1$; see \autoref{fig2}. \\\\
\textbf{Remark}. The quantity $M_1$ is well-defined since $c_g(\mu_0+m) \neq c_g(\mu_0+n)$. See the Appendix for the proof. The quantity $\gamma_0$ is also well-defined, as $\lambda_1 + i\mu_1 c_g(\mu_0+m)$ is guaranteed to be a complex number with nonzero real part. \\\\
Equating $\mu_1 =0$ maximizes the real part of $\lambda_1$ in \eqref{26a}. Thus, the Floquet exponent that corresponds to the most unstable spectral element of the high-frequency isola has asymptotic expansion
\begin{align} \label{27cc}
\mu_* = \mu_0 + \mathcal{O}(\varepsilon^2).
\end{align}
The corresponding real and imaginary components of this spectral element have asymptotic expansions
\begin{subequations} \label{27dd}
\begin{align}
\lambda_{*,r} &=  \varepsilon |\sigma| \sqrt{-(\mu_0+m)(\mu_0+n)} + \mathcal{O}(\varepsilon^2), \\
\lambda_{*,i} &= -\Omega(\mu_0+n) + \mathcal{O}(\varepsilon^2),
\end{align}
\end{subequations}
respectively. The former of these expansions provides an estimate for the growth rate of the $\Delta n = 1$ high-frequency instabilities. \autoref{fig2} compares these expansions with numerical results from FFH. Observe that, while the expansion for the real part is accurate, the expansion for the imaginary part requires a higher-order calculation; see Section 5.
 ~\\\\
If $\lambda$ is written as a sum of its real and imaginary components, $\lambda_{r}$ and $\lambda_{i}$, respectively, then eliminating dependence on $\mu_1$ between these quantities yields
\begin{align} \label{27d}
\frac{\lambda_{r}^2}{\varepsilon^2} +  \frac{\left(\lambda_{i}+\Omega(\mu_0+n)\right)^2}{\varepsilon^2\left( \frac{c_g(\mu_0+m)+c_g(\mu_0+n)}{c_g(\mu_0+m)-c_g(\mu_0+n)} \right)^2} = -\sigma^2 (\mu_0+m)(\mu_0+n) + \mathcal{O}(\varepsilon).
\end{align}
Thus, the $\Delta n = 1$ high-frequency isola is an ellipse to $\mathcal{O}(\varepsilon)$ with center at the collision site of eigenvalues $\lambda^{\mu_0}_{0,n}$ and $\lambda^{\mu_0}_{0,m}$ and with semi-major and -minor axes
\begin{subequations}
\begin{align}
a_1 &= \varepsilon|\sigma| \sqrt{-(\mu_0+m)(\mu_0+n)}, \\
b_1 &= a_1 \left| \frac{c_g(\mu_0+m)+c_g(\mu_0+n)}{c_g(\mu_0+m)-c_g(\mu_0+n)} \right|,
\end{align}
\end{subequations}
respectively. \\\\
Our asymptotic predictions agree well with numerics, particularly for the real component of the isola. There is some discrepancy between asymptotic and numerical results of the Floquet exponents and imaginary component of the isola, even when $\varepsilon = 10^{-3}$; see \autoref{fig3}. As noted before, this discrepancy is resolved in Section 5. \vspace*{0.4cm}
\subsection{High-Frequency Instabilities: $\Delta n = 2$ \nopunct} 
\begin{figure}[h!]
\centering \hspace*{-1cm} \includegraphics[width=15cm,height=7cm]{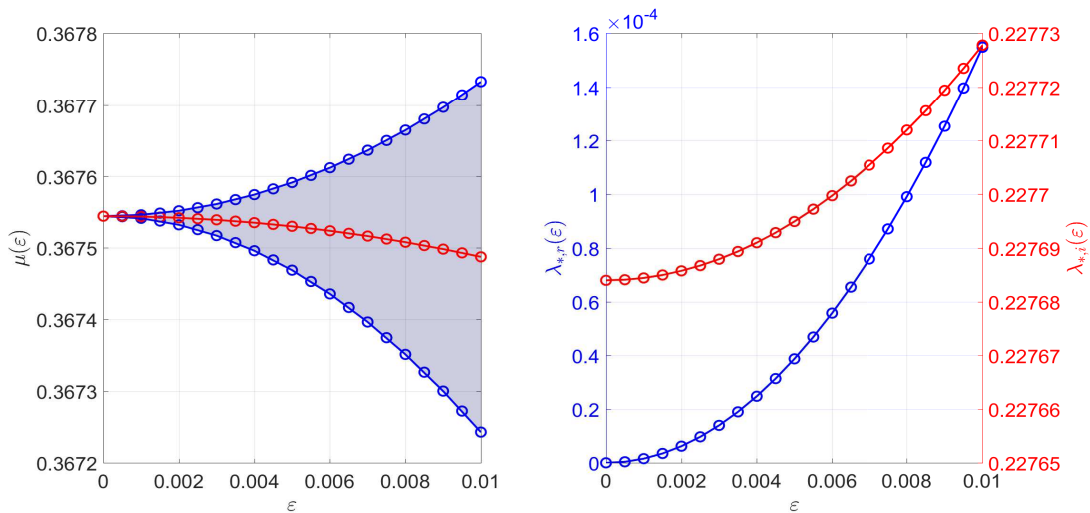}
\caption{(Left) Interval of Floquet exponents that parameterize the $\Delta n = 2$ high-frequency isola for parameters $\alpha = 1$, $\beta = 0.25$, and $\sigma = 1$ as a function of $\varepsilon$. ($\beta = 0.7$ only gives rise to a $\Delta n = 1$ isola: $\beta$ must be changed to satisfy \eqref{18} for a $\Delta n = 2$ isola to arise.) Solid blue curves indicate the asymptotic boundaries of this interval according to \eqref{35b}. Blue circles indicate the numerical boundaries computed using FFH. The solid red curve gives the Floquet exponent corresponding to the most unstable spectral element of the isola according to \eqref{35c}. The red circles indicate the same but computed numerically using FFH. (Right) The real (blue) and imaginary (red) parts of the most unstable spectral element of the isola as a function of $\varepsilon$. Solid curves illustrate asymptotic result \eqref{35d}. Circles illustrate results of FFH.}
\label{fig4}
\end{figure}

\begin{figure}[h!]
\centering \hspace*{-1cm} \includegraphics[width=15cm,height=7cm]{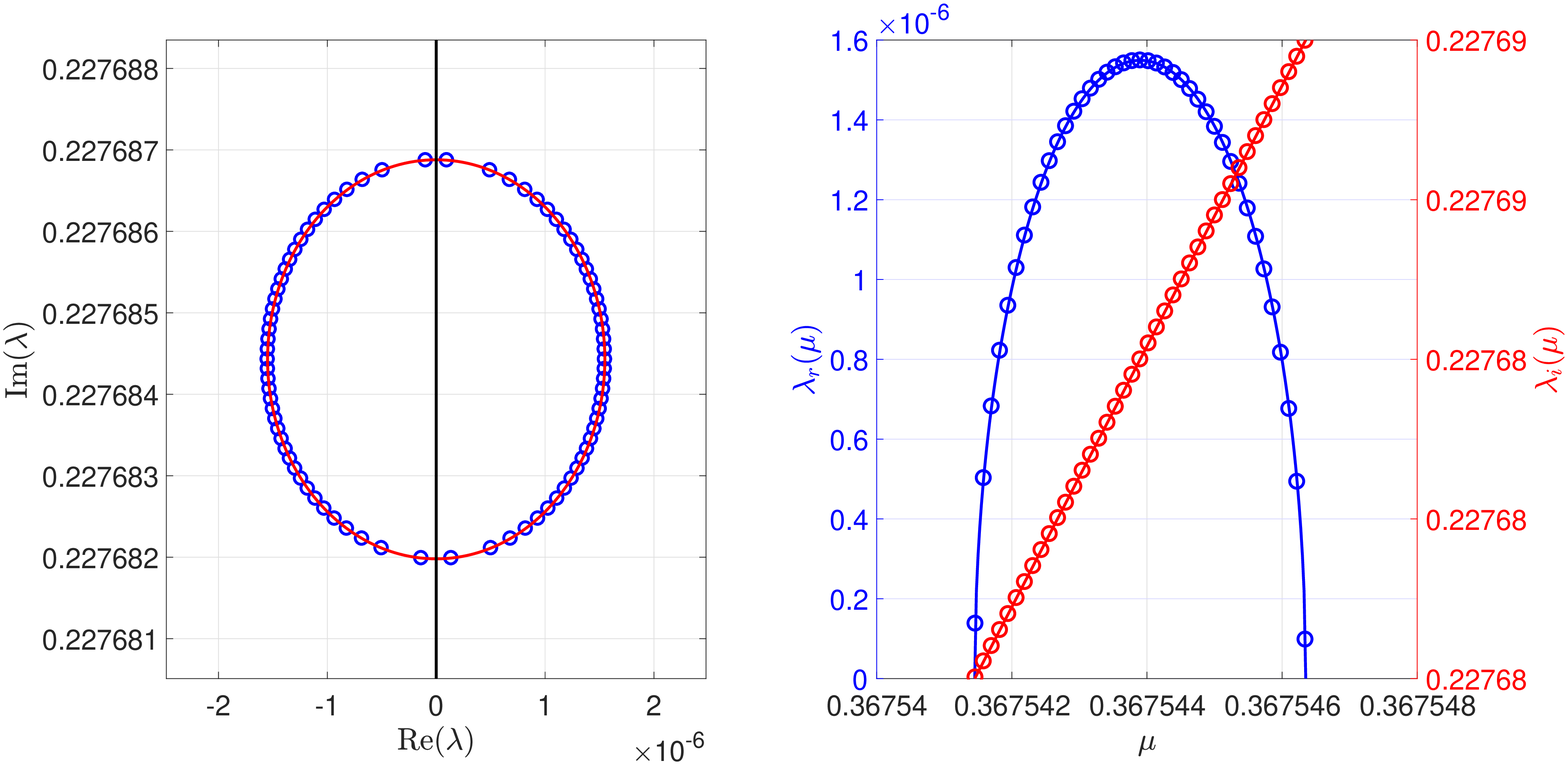}
\caption{(Left) $\Delta n = 2$ high-frequency isola for $\alpha = 1$, $\beta = 0.25$, $\sigma = 1$, and $\varepsilon = 10^{-3}$. The solid red curve is ellipse \eqref{36}. Blue circles are a subset of spectral elements from the numerically computed isola using FFH. Blue circles are a subset of spectral elements from the numerically computed isola using FFH. (Right) Floquet parameterization of the real (blue) and imaginary (red) parts of the isola. Solid curves illustrate asymptotic result \eqref{34a}. Circles indicate results of FFH. }
\label{fig5}
\end{figure}

~\\\\
\indent Suppose $m$, $n$, and $\mu_0$ satisfy the collision condition \eqref{17} for $\Delta n = 2$ and appropriately chosen $\beta$ parameter. Then, \eqref{20} gives a semi-simple eigenpair of $\mathcal{L}^{\mu_0}_0$, and we assume \eqref{21a}, \eqref{21b}, and \eqref{22} remain valid expansions for the eigenvalue, eigenfunction, and Floquet exponents in the vicinity of this semi-simple eigenpair, respectively. We obtain the coefficients of these expansions order by order in much the same way as for the $\Delta n = 1$ high-frequency instabilities. \\\\
\textit{The $\mathcal{O}(\varepsilon)$ Problem} \\\\
\indent  Substituting expansions \eqref{21a}, \eqref{21b}, and \eqref{22} into the spectral problem \eqref{15} and collecting $\mathcal{O}(\varepsilon)$ terms gives
\begin{align}
(\mathcal{L}^{\mu_0}_0-\lambda_0) w_1(x) =& -i\sigma (\mu_0+n-1)e^{i(n-1)x} +  \left[\lambda_1+i\mu_1 c_g(\mu_0+n)\right]e^{inx} \label{28} \\
&-i\sigma (\mu_0+n+1)(1+\gamma_0)e^{i(n+1)x} + \gamma_0\left[\lambda_1+i\mu_1c_g(\mu_0+m)\right]e^{imx}  \nonumber \\
&-i\sigma(\mu_0+m+1)e^{i(m+1)x}, \nonumber
\end{align}
where we have used \eqref{11} to replace $u_1(x)$. Though equation \eqref{28} shares similar features with \eqref{25}, $m - n \neq 1$ in this case. Thus, \eqref{28} cannot be simplified further. \\\\
The solvability conditions of \eqref{28} are
\begin{subequations}
\begin{align}
\lambda_1 + i\mu_1 c_g(\mu_0+n) &= 0, \\
\gamma_0 \left[\lambda_1 + i\mu_1 c_g(\mu_0+m) \right]   &= 0.
\end{align}
\end{subequations}
Since $c_g(\mu_0+m) \neq c_g(\mu_0+n)$ by the corollary provided in the Appendix and $\gamma_0 \neq 0$\footnote{Otherwise, our unperturbed eigenfunction $w_0(x)$ is not a superposition of two distinct Fourier modes.}, we must have
\begin{align} \label{29}
\lambda_1 = \mu_1 = 0.
\end{align}
Solving \eqref{28} for $w_1(x)$ by the the method of undetermined coefficients, one finds
\begin{align}  \label{30}
w_1(x) = \tau_{1,n-1}e^{i(n-1)x} + \tau_{1,n+1} e^{i(n+1)x} + \tau_{1,m+1} e^{i(m+1)x} + \gamma_1 e^{imx},
\end{align}
where $\gamma_1$ is a constant to be determined at higher order,
\begin{subequations}
\begin{align}
\tau_{1,n-1} &= Q_{n,n-1}, \\
\tau_{1,n+1} &= (1+\gamma_0) Q_{n,n+1}, \\
\tau_{1,m+1} &= \gamma_0Q_{n,m+1},
\end{align}
\end{subequations}
and
\begin{align} \label{29b}
Q_{N,M} = \frac{\sigma(\mu_0+M)}{\Omega(\mu_0+M)-\Omega(\mu_0+N)}.
\end{align}
Note that $w_1(x)$ does not have an $n^{\textrm{th}}$ Fourier mode, which is a consequence of the normalization \eqref{22b}. \\\\
\textit{The $\mathcal{O}(\varepsilon^2)$ Problem} \\\\
Substituting \eqref{21a}, \eqref{21b}, and \eqref{22} into \eqref{15} and collecting terms of $\mathcal{O}(\varepsilon^2)$ yields
\begin{align} \label{31}
(\mathcal{L}_0^{\mu_0}-\lambda_0)w_1(x) = \lambda_2w_0(x) - \mathcal{L}_2|_{\mu_1=0}w_0(x) - \mathcal{L}_1|_{\mu_1=0}w_1(x),
\end{align}
where $\mathcal{L}_1|_{\mu_1=0}$ is as before (but evaluated at $\mu_1=0$) and
\begin{align} \label{32b}
\mathcal{L}_2|_{\mu_1=0} = ic_0\mu_2 + c_2(i\mu_0+\partial_x) +3\mu_2 i(i\mu_0+\partial_x)^2   + 5\mu_2 i(i\mu_0+\partial_x)^4  + 2\sigma u_2(x)(i\mu_0+\partial_x) + 2\sigma u_2'(x).
\end{align}
As in the previous order, we evaluate the RHS of \eqref{31} using \eqref{11} to replace $u_1(x)$ and $u_2(x)$, \eqref{23} to replace $w_0(x)$, \eqref{30} to replace $w_1(x)$, and $m - n = 2$ to combine terms with exponential arguments proportional to $m-1$ and $n+1$. After some work, one arrives at the solvability conditions\footnote{To obtain \eqref{33a} and \eqref{33b}, one also needs evenness of $u_2(x)$ so that $\widehat{u_2(x)}_{-2} = \widehat{u_2(x)}_{2}$.}:
\begin{subequations}
\begin{align} \label{33a}
\lambda_2 + i\mathcal{C}_{\mu_2,\mu_0}^{n} &= i\gamma_0\mathcal{S}_2(\mu_0+n), \\
\gamma_0 \left[\lambda_2 + i\mathcal{C}_{\mu_2,\mu_0}^{m}\right] &= i\mathcal{S}_2(\mu_0+m), \label{33b}
\end{align}
\end{subequations}
where
\begin{subequations}
\begin{align}
\mathcal{C}_{\mu_2,\mu_0}^{N} &= \mu_2 c_g(\mu_0+N) - \mathcal{P}^{N}_{\mu_0}, \\
\mathcal{P}^N_{\mu_0} &= (\mu_0+N)\left[\sigma(Q_{n,N-1} + Q_{n,N+1} + 2\widehat{u_2(x)}_0) + c_2 \right], \\
\mathcal{S}_{2} &= \sigma(Q_{n,n+1}+2\widehat{u_2(x)}_2). \label{S2}
\end{align}
\end{subequations}
Similar to the $\Delta n = 1$ case, \eqref{33a} and \eqref{33b} are a nonlinear system for $\lambda_2$ and $\gamma_0$. The solution of this system is
\begin{subequations}
\begin{align}
\lambda_2 =& ~ -i\left(\frac{\mathcal{C}_{\mu_2,\mu_0}^m+\mathcal{C}_{\mu_2,\mu_0}^n}{2} \right) \pm \sqrt{-\left[\frac{\mathcal{C}_{\mu_2,\mu_0}^m-\mathcal{C}_{\mu_2,\mu_0}^n}{2} \right]^2 - \mathcal{S}_2^2(\mu_0+m)(\mu_0+n)},   \label{34a}\\
\gamma_0 =& ~ \frac{i\sigma(\mu_0+m)(Q_{n,n+1}+2\upsilon_{2,-2})}{\lambda_2 + i\mathcal{C}_{\mu_2,\mu_0}^{m}}. \label{34b}
\end{align}
\end{subequations}
Provided $\mathcal{S}_2 \neq 0$, there exists an interval of $\mu_2 \in (M_{2,-},M_{2,+})$, where
\begin{align} \label{35}
M_{2,\pm} = \frac{\mathcal{P}^m_{\mu_0} - \mathcal{P}^n_{\mu_0}}{c_g(\mu_0+m)-c_g(\mu_0+n)} \pm 2 \left| \frac{\mathcal{S}_2}{c_g(\mu_0+m)-c_g(\mu_0+n)} \right| \sqrt{-(\mu_0+m)(\mu_0+n)},
\end{align}
such that $\lambda_2$ has a nonzero real part. It is shown in the Appendix that $\mathcal{S}_2 \neq 0$ for all relevant values of $\beta$. Thus, the interval of Floquet exponents that parameterizes the $\Delta n = 2$ high-frequency isola to $\mathcal{O}(\varepsilon^2)$ is
\begin{align} \label{35b}
\mu \in \left(\mu_0 + \varepsilon^2 M_{2,-}, \mu_0 + \varepsilon^2 M_{2,+}\right).
\end{align}
Unlike when $\Delta n = 1$ \eqref{27b}, the center of this interval changes at the same rate as its width, and this width is an order of magnitude smaller than for the $\Delta n = 1$ instabilities. This explains why numerical detection of $\Delta n = 2$ instabilities presents a greater challenge than for $\Delta n =1$ instabilities; see \autoref{fig4}.
 ~\\\\
From the results above, we obtain an asymptotic expansion for the Floquet exponent of the most unstable spectral element of the $\Delta n = 2$ high-frequency isola:
\begin{align} \label{35c}
\mu_* = \mu_0 +  \frac{\mathcal{P}^m_{\mu_0} - \mathcal{P}^n_{\mu_0}}{c_g(\mu_0+m)-c_g(\mu_0+n)} \varepsilon^2+ \mathcal{O}(\varepsilon^3).
\end{align}
Asymptotic expansions for the real and imaginary component of this spectral element are
\begin{subequations} \label{35d}
\begin{align}
\lambda_{*,r} &=  \varepsilon^2 |\mathcal{S}_2| \sqrt{-(\mu_0+m)(\mu_0+n)} + \mathcal{O}(\varepsilon^3), \\
\lambda_{*,i} &=-\Omega(\mu_0+n)-\left[ \frac{\mathcal{P}_{\mu_0}^mc_g(\mu_0+n) - \mathcal{P}^{n}_{\mu_0}c_g(\mu_0+m) }{c_g(\mu_0+m)-c_g(\mu_0+n)}\right]\varepsilon^2  + \mathcal{O}(\varepsilon^3).
\end{align}
\end{subequations}
These expansions are in excellent agreement with numerical computations from FFH, as is seen in \autoref{fig4}. This is a consequence of resolving quadratic corrections to the real and imaginary components of $\Delta n = 2$ high-frequency isolas simultaneously, unlike in the $\Delta n = 1$ case.  \\\\
Analogous to the derivation of \eqref{27d}, the ellipse given by
\begin{align} \label{36}
\frac{\lambda_{r}^2}{\varepsilon^4}+& \frac{\left[\lambda_{i} +\Omega(\mu_0+n)+ \varepsilon^2 \left( \frac{\mathcal{P}_{\mu_0}^mc_g(\mu_0+n) - \mathcal{P}^{n}_{\mu_0}c_g(\mu_0+m) }{c_g(\mu_0+m)-c_g(\mu_0+n)}\right)\right]^2}{\varepsilon^4\left(\frac{c_g(\mu_0+m)+c_g(\mu_0+n)}{c_g(\mu_0+m)-c_g(\mu_0+n)}\right)^2 } =-\mathcal{S}_2^2(\mu_0+m)(\mu_0+n) + O(\varepsilon).
\end{align}
is asymptotic to the $\Delta n = 2$ high-frequency isola. This ellipse has center that drifts from the collision site at a rate comparable to its semi-major and -minor axes,
\begin{subequations}
\begin{align}
a_2 &= \varepsilon^2 |\mathcal{S}_2| \sqrt{-(\mu_0+m)(\mu_0+n)} \\
b_2 &= a_2 \left| \frac{c_g(\mu_0+m)+c_g(\mu_0+n)}{c_g(\mu_0+m)-c_g(\mu_0+n)} \right|,
\end{align}
\end{subequations}
respectively. This behavior contrasts with that seen in the $\Delta n = 1$ case, where the center drifts slower than the axes grow. Comparison with numerical computations using FFH show that \eqref{36} is an excellent approximation for $\Delta n = 2$ high-frequency isolas; see \autoref{fig5}.  ~\\\\
\subsection{High-Frequency Instabilities: $\Delta n \geq 3$ \nopunct} ~\\\\
The approach used to obtain leading-order behavior of the $\Delta n = 1,2$ high-frequency isolas generalizes to higher-order isolas. The method consists of the following steps, each of which is readily implemented in a symbolic programming language: ~\\
\begin{enumerate}[label=(\roman*)]
\item  Given $\Delta n \in \mathbb{N}$, determine the unique $\mu_0$, $m$, and $n$ to satisfy collision condition \eqref{17}, assuming $\beta$ satisfies \eqref{18}.
\item  Expand about the collided eigenvalues in a formal power series of $\varepsilon$ and similarly expand their corresponding eigenfunctions and Floquet exponents. To maintain uniqueness of the eigenfunctions, choose the normalization \eqref{22b}.
\item Substitute these expansions into the spectral problem \eqref{15}. Collecting like powers of $\varepsilon$, construct a hierarchy of inhomogeneous linear problems to solve.
\item Proceed order by order. At each order, impose solvability and normalization conditions. Invert the linear operator against its range using the  method of undetermined coefficients. Use previous normalization and solvability conditions as well as the collision condition to simplify problems if necessary.
\end{enumerate} ~\\
We conjecture that this method yields the first nonzero real part correction to the $\Delta n^{\textrm{th}}$ high-frequency isola at $\mathcal{O}(\varepsilon^{\Delta n})$. We have shown that this conjecture holds for $\Delta n = 1,2$. For $\Delta n = 3$, one can show that the high-frequency isola is asymptotic to the ellipse
\begin{align} \label{37}
\frac{\lambda_r^2}{\varepsilon^6} + \frac{\left(\lambda_i + \Omega(\mu_0+n) +\varepsilon^2\left[\frac{\mathcal{P}^{m}_{\mu_0}c_g(\mu_0+n) - \mathcal{P}^{n}_{\mu_0} c_g(\mu_0+m)}{c_g(\mu_0+m)-c_g(\mu_0+n)} \right]     \right)^2}{\varepsilon^6 \left(\frac{c_g(\mu_0+m)+c_g(\mu_0+n)}{c_g(\mu_0+m)-c_g(\mu_0+n)} \right)^2} = -\mathcal{S}_3^2(\mu_0+m)(\mu_0+n) + O(\varepsilon),
\end{align}
where
\begin{align}
\mathcal{S}_3 = \sigma \left[ Q_{n,n+1} Q_{n,n+2} + 2\widehat{u_2(x)}_{2}(Q_{n,n+1} + Q_{n,n+2}) + 2\widehat{u_3(x)}_{3} \right].
\end{align}
The semi-major and -minor axes  of \eqref{37} scale as $\mathcal{O}(\varepsilon^3)$, as the conjecture predicts.  If true for all $\Delta n$, this conjecture explains why higher-order isolas are challenging to detect both in numerical and perturbation computations of the stability spectrum.  \\\\
One notices that the center of \eqref{37} drifts similarly to that of the $\Delta n = 2$ high-frequency isolas \eqref{36}. In fact, centers of higher-order isolas (beyond $\Delta n = 1$) all drift at a similar rate, as these isolas all satisfy the same $\mathcal{O}(\varepsilon^2)$ problem and, hence, yield corrections at this order. Consequently, one can expect to incur corrections to the imaginary component of the high-frequency isola before reaching $\mathcal{O}(\varepsilon^{\Delta n})$, making it more difficult to prove our conjecture about the first occurence of a nonzero real part correction.

\section{$\Delta n = 1$ High-Frequency Instabilities at Higher-Order}

As we saw in Section 4.1, the asymptotic formulas for the $\Delta n = 1$ high-frequency isola fail to capture its $\mathcal{O}(\varepsilon^2)$ drift along the imaginary axis. This is expected, as we only considered the $\mathcal{O}(\varepsilon)$ problem. In this section, we go beyond the leading-order behavoir of these instabilities. We expect similar calculations to arise if one considered the $\mathcal{O}(\varepsilon^p)$ problem for a generic $\Delta n$ isola, where $p > \Delta n$.  \\
\subsection{The $\mathcal{O}(\varepsilon)$ Problem Revisited \nopunct} ~\\\\
Finishing our work from Section 4.1, we solve for $w_1(x)$ in \eqref{25}. We find
\begin{align}
w_1(x) = Q_{n,n-1} e^{i(n-1)x} + \gamma_0 Q_{n,m+1} e^{i(m+1)x} + \gamma_1 e^{imx},
\end{align}
where $Q_{N,M}$ is defined as in \eqref{29b} and $\gamma_1$ is an undetermined constant at this order. \\
\subsection{The $\mathcal{O}(\varepsilon^2)$ Problem \nopunct} ~\\\\At $\mathcal{O}(\varepsilon^2)$, we have
\begin{align} \label{40}
(L_0-\lambda_0)w_2(x) =&~\lambda_2 w_0(x) + \lambda_1 w_1(x) - i c_0(\mu_1w_1(x) + \mu_2 w_0(x)) -c_2(i\mu_0+\partial_x)w_0(x) \\
&-3i(i\mu_0+\partial_x)^2(\mu_1w_1(x) + \mu_2 w_0(x)) + 3\mu_1^2(i\mu_0+\partial_x)w_0(x) \nonumber \\ &- 5\beta i(i\mu_0+\partial_x)^4(\mu_1w_1(x)+\mu_2w_0(x))+ 10\beta\mu_1^2(i\mu_0+\partial_x)^3w_0(x) \nonumber \\ &- 2\sigma (i\mu_0 + \partial_x)(u_{1}(x)w_1(x) + u_{2}(x)w_0(x)) - 2\sigma i \mu_1 u_{1}(x)w_0(x).
\nonumber\end{align}
After substituting $w_0(x)$, $w_1(x)$, $u_{1}(x)$, and $u_{2}(x)$ into \eqref{40}, solvability conditions of the second-order problem form the linear system
\begin{align} \displaystyle
\begin{pmatrix} 1 & -i\sigma(\mu_0+n) \\ \gamma_0 & \lambda_1 + i\mu_1c_g(\mu_0+m) \end{pmatrix} \begin{pmatrix} \lambda_2 \\ \gamma_1 \end{pmatrix} = i \begin{pmatrix} \sigma \gamma_0 \mu_1 - \tilde{\mathcal{C}}_{\mu_2,\mu_1,\mu_0}^{n,-1}   \\ \sigma\mu_1-\gamma_0 \tilde{\mathcal{C}}_{\mu_2,\mu_1,\mu_0}^{m,1}  \end{pmatrix}, \label{40.5}
\end{align}
where
\begin{subequations}
\begin{align}
\tilde{C}_{\mu_2,\mu_1,\mu_0}^{N,k} &= \mu_2 c_g(\mu_0+N) - \tilde{\mathcal{P}}_{\mu_0}^{N,k} + \mu_1^2 \mathcal{D}^N_{\mu_0}, \\
\tilde{P}_{\mu_0}^{N,k} &= (\mu_0+N)\left[\sigma(Q_{n,N+k} +2\upsilon_{2,0})+c_2 \right], \\
D_{\mu_0}^{N} &= 3(\mu_0+N)-10\beta(\mu_0+N)^3.
\end{align}
\end{subequations}
For $\mu_1 \in (-M_1,M_1)$ \eqref{27a}, one can show that
\begin{align}
\textrm{det}\begin{pmatrix} 1 & -i\sigma(\mu_0+n) \\ \gamma_0 & \lambda_1 + i\mu_1 c_g(\mu_0+m) \end{pmatrix}  &= 2\lambda_{1,r}.
\end{align}
Since $\lambda_{1,r} \neq 0$ for this interval of $\mu_1$, it follows that \eqref{40.5} is invertible. Using Cramer's rule and \eqref{26a} and \eqref{26b}, the solvability conditions at $\mathcal{O}(\varepsilon)$, gives
\begin{align} \label{44}
\lambda_{2} = -\frac{i}{2\lambda_{1,r}} \left(\mathcal{A}\lambda_1 +i\mu_1 \mathcal{B} \right),
\end{align}
where
\begin{subequations}
\begin{align}
\mathcal{A} &= \tilde{\mathcal{C}}_{\mu_2,\mu_1,\mu_0}^{m,1} + \tilde{\mathcal{C}}_{\mu_2,\mu_1,\mu_0}^{n,-1}, \\
\mathcal{B} &= c_g(\mu_0+m) \tilde{\mathcal{C}}_{\mu_2,\mu_1,\mu_0}^{n,-1} + c_g(\mu_0+n)\tilde{\mathcal{C}}_{\mu_2,\mu_1,\mu_0}^{m,1}- \sigma^2(2\mu_0+m+n). \\ \nonumber
\end{align}
\end{subequations}
\subsection{Determination of $\mu_2$: The Regular Curve Condition \nopunct} ~\\\\
A quick calculation shows that $\lambda_2$ has two branches, $\lambda_{2,+}$ and $\lambda_{2,-}$, and, for any $\mu_2 \in \mathbb{R}$, $\lambda_{2,+} = -\overline{\lambda_{2,-}}$. Consequently, \eqref{44} results in a spectrum that is symmetric about the imaginary axis regardless of $\mu_2$. We want this spectrum to be a continuous, closed curve about the imaginary axis. As we shall see, this additional constraint is enough to determine $\mu_2$ uniquely. We call this additional constraint the \textit{regular curve condition}. \\\\
To motivate the regular curve condition, consider the real and imaginary parts of \eqref{44}:
\begin{subequations}
\begin{align}
\lambda_{2,r} &= \frac{1}{2\lambda_{1,r}} \left(\mathcal{A}\lambda_{1,i} + \mu_1 \mathcal{B} \right), \label{44b} \\
\lambda_{2,i} &= -\frac{\mathcal{A}}{2}.
\end{align}
\end{subequations}
As $|\mu_1|$ approaches $M_1$, $\lambda_{1,r}$ approaches zero. To avoid unwanted blow-up of $\lambda_{2,r}$, we must impose
\begin{align} \label{45}
\displaystyle \lim_{|\mu_1| \rightarrow M_1} \left(\mathcal{A}\lambda_{1,i} + \mu_1 \mathcal{B} \right) = 0.
\end{align}
Since $\mu_1$ appears in $\mathcal{A}$ only as $\mu_1^2$, we can rewrite \eqref{45} with the help of \eqref{26c} as
\begin{align} \label{48}
 \lim_{\mu_1^2 \rightarrow M_1^2} \left(-\frac{\mathcal{A}}{2}(c_g(\mu_0+m)+c_g(\mu_0+n)) + \mathcal{B} \right) = 0.
\end{align}
Equation \eqref{48} is the \textit{regular curve condition} for second-order corrections to the $\Delta n = 1$ isola. Unpacking the definitions of $\mathcal{A}$ and $\mathcal{B}$ above, the \textit{regular curve condition} implies that
\begin{align} \label{49}
\mu_2 = \frac{\textrm{P}^{m,1}_{\mu_0}-\textrm{P}^{n,-1}_{\mu_0}}{c_g(\mu_0+m)-c_g(\mu_0+n)} - \frac{2\sigma^2(2\mu_0+m+n)}{(c_g(\mu_0+m)-c_g(\mu_0+n))^2},
\end{align}
where
\begin{align}
\textrm{P}^{N,k}_{\mu_0} = \tilde{\mathcal{P}}^{N,k}_{\mu_0} + \frac{4\sigma^2(\mu_0+m)(\mu_0+n)}{(c_g(\mu_0+m)-c_g(\mu_0+n))^2} \mathcal{D}^{N}_{\mu_0}.
\end{align}
Therefore, to $\mathcal{O}(\varepsilon^2)$, the asymptotic interval of Floquet exponents that parameterizes the $\Delta n = 1$ high-frequency isola is
\begin{align} \label{49b}
\mu \in \left(\mu_0 - \varepsilon M_1 + \varepsilon^2 \mu_2, \mu_0 + \varepsilon M_1 + \varepsilon^2 \mu_2\right).
\end{align} ~\\
\subsection{The Most Unstable Eigenvalue \nopunct} 
\begin{figure}[t!]
\centering \hspace*{-1cm} \includegraphics[width=15cm,height=7cm]{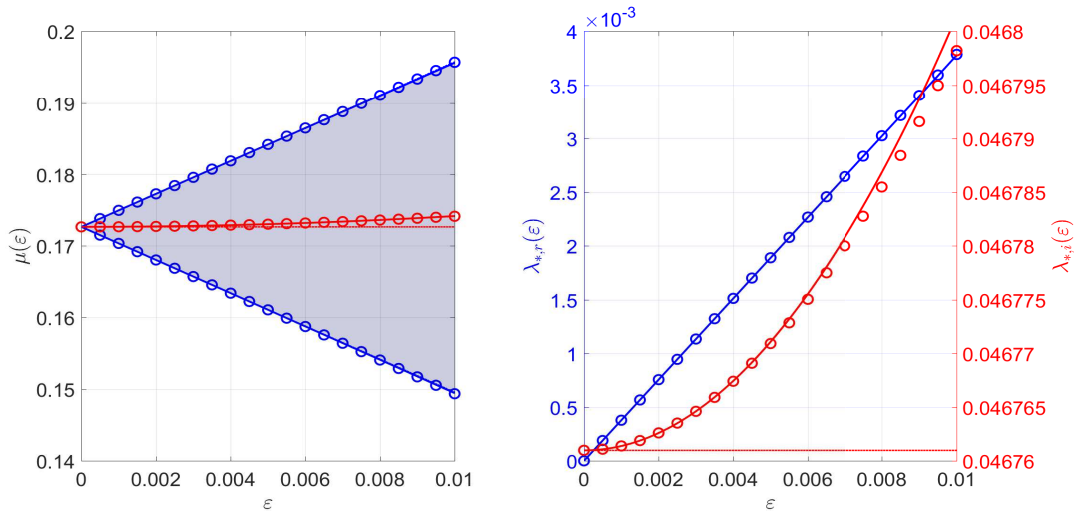}
\caption{(Left) Interval of Floquet exponents that parameterize the $\Delta n = 1$ high-frequency isola for parameters $\alpha = 1$, $\beta = 0.7$, and $\sigma = 1$ as a function of $\varepsilon$. Solid blue curves indicate the asymptotic boundaries of this interval according to \eqref{49b}, while the dotted blue curves give the $\mathcal{O}(\varepsilon)$ result. Blue circles indicate the numerical boundaries computed using FFH. The solid red curve gives the Floquet exponent corresponding to the most unstable spectral element of the isola according to \eqref{56}, while the dotted red gives the $\mathcal{O}(\varepsilon)$ result. The red circles indicate the same but computed numerically using FFH. (Right) The real (blue) and imaginary (red) parts of the most unstable spectral element of the isola as a function of $\varepsilon$. Solid curves illustrate asymptotic result \eqref{56b}. Dotted curves illustrate the asymptotic results only to $\mathcal{O}(\varepsilon)$. Circles illustrate results of FFH.}
\label{fig6}
\end{figure}

\begin{figure}[t!]
\centering \hspace*{-1cm} \includegraphics[width=15cm,height=7cm]{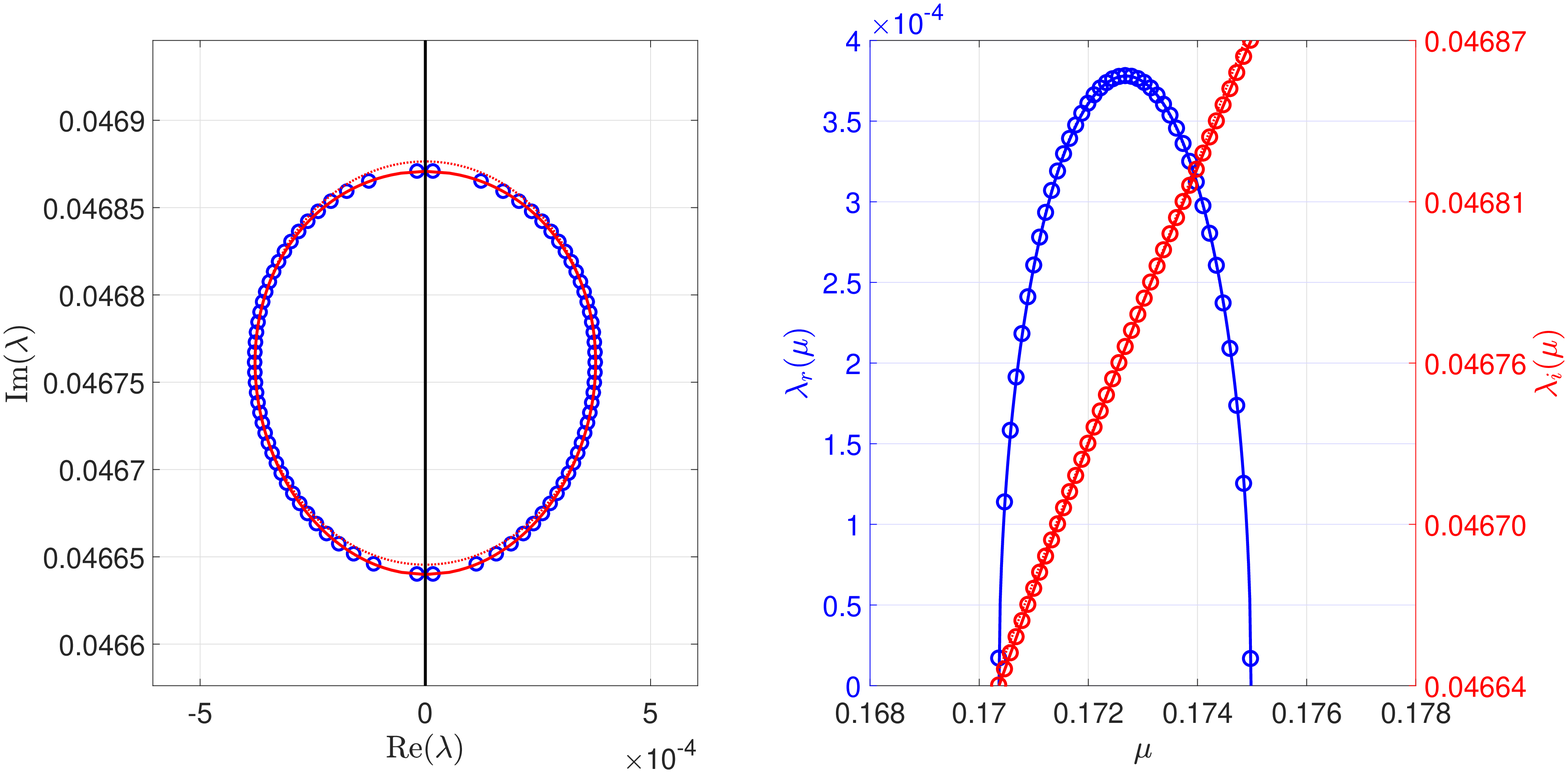}
\caption{(Left) $\Delta n = 1$ high-frequency isola for $\alpha = 1$, $\beta = 0.7$, $\sigma = 1$, and $\varepsilon = 10^{-3}$. The solid red curve is parameterized by \eqref{44}. This curve is no longer an ellipse, but a more complicated algebraic curve. For comparison, the dotted red curve is the ellipse found at $\mathcal{O}(\varepsilon)$. Blue circles are a subset of spectral elements from the numerically computed isola using FFH. (Right) Floquet parameterization of the real (blue) and imaginary (red) parts of the isola. Solid curves illustrate asymptotic result \eqref{44}. Dotted curves illustrate the asymptotic results only to $\mathcal{O}(\varepsilon)$. Circles indicate results of FFH. }
\label{fig7}
\end{figure}

~\\\\
To $\mathcal{O}(\varepsilon)$, the expression for the real part of the $\Delta n = 1$ high-frequency isola is
\begin{align} \label{50}
\lambda^{(1)}_{r} := \varepsilon \lambda_{1,r} =  \pm \varepsilon \sqrt{-\mu_1^2\left[\frac{c_g(\mu_0+m)-c_g(\mu_0+n)}{2} \right]^2-\sigma^2(\mu_0+m)(\mu_0+n)}.
\end{align}
The most unstable eigenvalue of \eqref{50} occurs when $\mu_1 = \mu_{*,1}$, where  $\mu_{*,1}$ is a critical point of $\lambda^{(1)}_{r}$:
\begin{align} \label{51}
\frac{\partial \lambda^{(1)}_{r}}{\partial \mu_1} \biggr|_{\mu_{*,1}} = 0.
\end{align}
Solving \eqref{51}, one finds $\mu_{*,1} = 0$, and we conclude that the Floquet exponent that corresponds to the most unstable eigenvalue is $\mu_* = \mu_0 + O(\varepsilon^2)$, as found in Section 4.1. \\\\
To $\mathcal{O}(\varepsilon^2)$, the real part of our isola is
\begin{equation} \label{53}
\lambda^{(2)}_r := \varepsilon \lambda_{1,r} + \varepsilon^2 \lambda_{2,r},
\end{equation}
where $\lambda_{1,r}$ is given in \eqref{50} and $\lambda_{2,r}$ is given in \eqref{44b}. Without loss of generality, we choose the positive branch of $\lambda_{1,r}$. \\\\  Taking inspiration from \eqref{51}, we consider the critical points of \eqref{53}:
\begin{equation}
\label{54}
\frac{\partial \lambda_{r}^{(2)}}{\partial \mu_1} \biggr|_{\mu_{*,1}} = 0.
\end{equation}
After some tedious calculations, \eqref{54} yields the following equation for $\mu_{*,1}$:
\begin{equation} \label{55}
\begin{aligned}
-\mu_{*,1} \lambda_{*,1,r}^2\left(\frac{c_g(\mu_0+m)-c_g(\mu_0+n)}{2} \right)^2+ \frac{\varepsilon}{2} \biggr[ \lambda_{*,1,r}^{2} \biggr( \lambda_{*,1,i} \frac{\partial \mathcal{A}}{\partial \mu_1} \biggr|_{\mu_{*,1}} + \mathcal{A}_* \frac{\partial \lambda_{1,i}}{\partial \mu_1 } \biggr|_{\mu_{*,1}} \phantom) \phantom] \\ \phantom[ \phantom( + \mu_{*,1} \frac{\partial \mathcal{B}}{\partial \mu_1} \biggr|_{\mu_{*,1}} + \mathcal{B}_*\biggr) + \mu_{*,1} \left( \mathcal{A}_* \lambda_{*,1,i} + \mu_{*,1} \mathcal{B}_* \right) \left(\frac{c_g(\mu_0+m) - c_g(\mu_0+n)}{2}\right)^2 \biggr] = 0, 
\end{aligned}
\end{equation}
where it is understood that starred variables are evaluated at $\mu_{*,1}$. Unpacking the definitions of $\mathcal{A}$, $\mathcal{B}$, $\lambda_{1,r}$, and $\lambda_{1,i}$ reveals that \eqref{55} is a quartic equation for $\mu_{*,1}$ with the highest degree coefficient multiplied by the small parameter $\varepsilon$. Rather than solving for $\mu_{*,1}$ directly, we obtain the roots perturbatively. \\\\
An application of the method of dominant balance to \eqref{55} shows that all of its roots have leading order behavior $\mathcal{O}(\varepsilon^{-1})$, except for one. Because we anticipate that $\lim_{\varepsilon \rightarrow 0} \mu_{*,1} = 0$ to match results at the previous order,  it is this non-singular root that we expect to yield the next order correction for $\mu_{*,1}$. Therefore, we need not concern ourselves with singular perturbation methods and, instead, make the following ansatz:
\begin{align}
\mu_{*,1} = 0 + \varepsilon \mu_{*,1,1} + O(\varepsilon^2).
\end{align}
Plugging our ansatz into \eqref{55} and keeping terms of lowest power in $\varepsilon$, we find the following linear equation to solve for $\mu_{*,1,1}$:
\begin{align}
-\mu_{*,1,1} \left( \frac{c_g(\mu_0+m)-c_g(\mu_0+n)}{2} \right)^2 + \frac12 \left( \mathcal{B}_0 - \mathcal{A}_0 \right) = 0,
\end{align}
where $\mathcal{A}_0$ and $\mathcal{B}_0$ are $\mathcal{A}$ and $\mathcal{B}$ evaluated at $\mu_1 = 0$, respectively. Using the definition of $\mathcal{A}$ and $\mathcal{B}$ together with the expression for $\mu_2$ in \eqref{49} above, one finds that
\begin{align}
\mu_{*,1,1} = -4\sigma^2(\mu_0+m)(\mu_0+n) \left(\frac{\mathcal{D}^m_{\mu_0} - \mathcal{D}^n_{\mu_0}}{(c_g(\mu_0+m)-c_g(\mu_0+n))^3}\right).
\end{align}
It follows that the Floquet exponent corresponding to the most unstable eigenvalue of the $\Delta n = 1$ high-frequency isola is
\begin{align} \label{56}
\mu_* = \mu_0 + \varepsilon^2(\mu_2 + \mu_{*,1,1})+O(\varepsilon^3).
\end{align}
The most unstable eigenvalue is then
\begin{align} \label{56b}
\lambda_* = \lambda_0 + \varepsilon \lambda_1|_{\mu_1 = 0}+ \varepsilon^2 \lambda_2|_{\mu_1 = 0, \mu_2 =\mu_2+\mu_{*,1,1}} + O(\varepsilon^3).
\end{align}
\autoref{fig6} and \autoref{fig7} show improvements to results in \autoref{fig2} and \autoref{fig3}, respectively, as a result of our higher-order calculations. 
\section{Conclusions}
In this work, we investigate the asymptotic behavior of high-frequency instabilities of small-amplitude Stokes waves of the Kawahara equation. For the largest of these instabilities ($\Delta n = 1,2$), we introduce a perturbation method to compute explicitly
\begin{enumerate}[label=(\roman*)]
\item the asymptotic interval of Floquet exponents that parameterize the high-frequency isola,
\item the leading-order behavior of its most unstable spectral elements, and
\item the leading-order curve asymptotic to the isola.
\end{enumerate}
We outline the procedure to compute these quantities for higher-order isolas. For the first time, we compare these asymptotic results with numerical results and find excellent agreement between the two.  We also obtain higher-order asymptotic results for the $\Delta n = 1$ high-frequency isolas by introducing the \textit{regular curve condition}.  \\\\
The perturbation method used throughout our investigation holds only for nonresonant Stokes waves \eqref{10}. Resonant waves require a modified Stokes expansion, and as a result of this modification, the leading-order behavior of the high-frequency isolas will change. Some numerical work has been done to investigate this effect \cite{trichtchenko18}, but no perturbation methods have been proposed.  
\section{Appendix}
\noindent \textbf{Theorem 1}. For each $\Delta n \in \mathbb{N}$, if $\beta$ satisfies \eqref{10} and \eqref{18}, there exists a unique $\mu_0 \in [0,1/2]$ and unique $m\neq n \in \mathbb{Z}$ such that the collision condition \eqref{17} is satisfied. \\\\
\textit{Proof}.
Define \begin{align} F(k;\Delta n) = \frac{\Omega(k + \Delta n)- \Omega(k)}{\Delta n}. \end{align} Using the definition of the dispersion relation $\Omega$,
\begin{equation}
 \begin{aligned} F(k;\Delta n) =&~ 5\beta k^4 + 10\beta \Delta n k^3 + (10 \beta (\Delta n)^2 -3)k^2 +(5\beta (\Delta n)^3 - 3\Delta n)k \\&+ \beta (\Delta n)^4 - (\Delta n)^2 + 1 - \beta.  \end{aligned} \end{equation}
A direct calculation shows that
\begin{align}
F(k ;\Delta n) = F(-(k+\Delta n);\Delta n).
\end{align} Hence, the graph of $F$ is symmetric about $k = -\Delta n/2$. We prove the desired result for the various cases of $\Delta n$.  \\\\
\textbf{Case 1.} Suppose $\Delta n = 1$. Then, $k_1=0$ and $k_2=-1$ are roots of $F$ by inspection. The remaining roots are
\begin{align}
k_{3,4} = \frac{-1 \pm \sqrt{\frac{12}{5\beta}-3}}{2}.
\end{align}
Because $\beta$ satisfies \eqref{18}, one can show that \begin{align} 0<\frac{12}{5\beta}-3 < 1,\end{align} so that $k_{3,4} \in (k_2,k_1)$. Because $F$ is symmetric about $k = -\Delta n/2$, we have $k_3 \in \left(-1/2,0\right)$ and $k_4 \in \left(-1,-1/2\right)$. \\\\
Each of these wavenumbers $k_j$ is mapped to a Floquet exponent $\mu \in (-1/2, 1/2]$ according to
\begin{align}
\mu(k) = k - [k],
\end{align}
where $[\cdot]$ denotes the nearest integer function\footnote{Our convention is $\left[p/2\right] = (p-1)/2$ for $p$ odd.}. Both $k_1$ and $k_2$ map to $\mu=0$. One checks that $\mu(k_3)=-\mu(k_4)\neq 0$ and $|\mu(k_3)| = |\mu(k_4)|<1/2$, since $k_4 = -(k_3+1)$ and $-1/2 < k_3 < 0$. Thus, the requisite $\mu_0 \in (0,1/2)$ is $\mu(k_j)$, where $j$ is either $3$ or $4$ depending on which has the correct sign. Then, $n = [k_j]$ and $m = n + \Delta n$. These are unique by the uniqueness of $k_j$.  \\\\
\textbf{Case 2.} Suppose $\Delta n = 2$. A calculation of $F(-1;2)$ and $F_k(-1;2)$ shows that $k_{1,2} = -1$ is a double root. The remaining roots are
\begin{align}
k_{3,4} = -1 \pm \sqrt{\frac{3}{5\beta}-2}.
\end{align}
Clearly, $\mu(k_{1,2})=0$. Since $k_4 = -(k_3 + 2)$, we again have $\mu(k_3) = -\mu(k_4)$. Also, from the formula for $k_3$ above, we have that $-1 < k_3 < 0$ by \eqref{18}, so $\mu(k_3)$ is nonzero. Thus, $\mu(k_j)$ is the requisite $\mu_0 \in (0,1/2]$, where $j$ is either $3$ or $4$ depending on which has the correct sign. Again, $n = [k_j]$ and $m = n + \Delta n$ are uniquely defined. \\\\
\textbf{Remark}. Unlike in the first case, we cannot guarantee $\mu_0 \neq 1/2$. Indeed, this value can be achieved when $\beta = 4/15$. \\\\
\textbf{Case 3.} Suppose $\Delta n \geq 3$. The discriminant of $F(k;\Delta n)$ with respect to $k$ is
\begin{align}
\Delta_k[F] = 5\beta\left[(\Delta n)^2-4\right]\left[\beta((\Delta n)^2+4)-4\right]\left[5\beta\left(\beta\left((\Delta n)^4+4\right)-2\left((\Delta n)^2+2\right)\right)+9\right]^2.
\end{align}
For $\beta$ satisfying inequality \eqref{18}, we have $\Delta_k[F]<0$, implying there are two distinct real roots of $F$. These roots must be nonpositive by an application of Descartes' Rule of Signs on $F$. Without loss of generality, suppose $k_2 < k_1$. Then, by the symmetry of $F$ about $k = -\Delta n/2$, $k_2 = -(k_1+\Delta n)$. It follows that $\mu(k_1) =-\mu(k_2)$. Thus, $\mu(k_j)$ is the requisite value of $\mu_0 \in [0,\frac12 ]$, where $j = 1$ or $2$ depending on which has the correct sign. The integers $n$ and $m$ are uniquely defined as before.
\\\\
\textbf{Remark}. In \cite{kollar19}, $\beta = 1/((\Delta n/2)^2+1)$ is included in inequality \eqref{18} when $\Delta n \geq 3$. For this $\beta$, $F(k;\Delta n)$ has a double root at $k_* = -\Delta n/2$. However, $\Omega(k_*) = 0$, which corresponds to an eigenvalue collision at the origin. Eigenvalue collisions at the origin do not satisfy \eqref{17}.  \\\\
In each of these cases, we have found $k_j<0$ such that $F(k_j;\Delta n) = 0$. Importantly, one must check that $\Omega(k_j) \neq 0$ for such $k_j$. Indeed, suppose $\Omega(k_j) = 0$. A direct calculation shows that $k_j = \pm 1$, $0$, or $k_j^2 = (1-\beta)/\beta$. Clearly $k_j =0$ or $1$ contradict that $k_j<0$. If $k_j = - 1$, then $F(-1;\Delta n) = 0$ implies $\beta = 1/[(\Delta n-1)^2+1]$, which contradicts \eqref{10} when $\Delta n \neq 2$. If $\Delta n = 2$, $\beta = 1/2$, which contradicts \eqref{18}. \\\\
It remains to be seen if $k_j^2 = (1-\beta)/\beta$ leads to contradiction. Indeed, a straightforward (although tedious) calculation shows that, if $k_j^2 = (1-\beta)/\beta$, $F(k_j;\Delta n) = 0$ implies $\beta = 0$, $\beta = 1/[1+(\Delta n/2)^2]$, $\beta = 1/[1+(\Delta n - 1)^2]$, or $\beta = 1/[1+(\Delta n + 1)^2]$. All of these lead to contradictions of \eqref{10} or \eqref{18}. Therefore, we must have $\Omega(k_j) = \Omega(\mu_0+n) = \Omega(\mu_0+m) \neq 0$ in all cases, as desired.    \begin{flushright}$\square$ \end{flushright}
In expressions for the isolas derived in Sections 4 and 5, factors of $c_g(\mu_0+m) - c_g(\mu_0+n)$ appear in denominators. A consequence of Theorem 1 is that this factor is never zero:   \\\\
\textbf{Corollary 1.} Fix $\Delta n \in \mathbb{N}$ and choose $\beta$ to satisfy \eqref{10} and \eqref{18}. Consider $\mu_0 \in [0,1/2]$ that solves $\Omega(\mu_0+m) = \Omega(\mu_0+n)$ for unique integers $m,n$ such that $m = n + \Delta n$. Suppose, in addition, that $\mu_0$ solves $c_g(\mu_0+m) = c_g(\mu_0+n)$, where $c_g(k) = \Omega'(k)$.  Then, $\Delta n = 2$ and $\mu_0 = 0$. \\\\
\textit{Proof}.
If  $\Omega(\mu_0+m) = \Omega(\mu_0+n)$ and $c_g(\mu_0+m) = c_g(\mu_0+n)$, then $k_0 = \mu_0+m$ is a double root of $F(k;\Delta n )$. From the proof of the theorem above, the only double root is $k_0 = -1$ (i.e. $\mu_0 = 0$) when $\Delta n = 2$. \begin{flushright} $\square$ \end{flushright}
The corresponding eigenvalue collision for this $\Delta n$ and $\mu_0$ happens at the origin in the complex spectral plane and is not of interest to us. Thus, $c_g(\mu_0+m) \neq c_g(\mu_0+n)$. \\\\
In Section 4.2, the quantity $\mathcal{S}_2$ \eqref{S2} must be nonzero in order for $\Delta n = 2$ high-frequency instabilities to exist at $\mathcal{O}\left(\varepsilon^2\right)$. The following corollary shows $\mathcal{S}_2 \neq 0$ for $\beta$ satisfying inequality \eqref{18}. \\\\
\textbf{Corollary 2.} For $\mathcal{S}_2$ defined in \eqref{S2} and $\beta$ satisfying inequality \eqref{18} for $\Delta n = 2$, $\mathcal{S}_2 \neq 0$. \\\\
\textit{Proof}. Since $\Delta n = 2$, we have from \eqref{18} that $1/5 < \beta < 3/10$. In addition, from Theorem 1 and Corollary 1, we know $k_{1,2}=1$ is a double root of $F(k;\Delta n)$ for all $\beta$ in this interval, and the remaining roots of $F$ are
\begin{align}
k_{3,4} = -1 \pm \sqrt{\frac{3}{5\beta}-2}.
\end{align}
These remaining roots correspond to the nonzero eigenvalue collisions that give rise to the $\Delta n = 2$ high-frequency instability. \\\\
The quantity $\mathcal{S}_2$ can be written in terms of $k_{3,4}$ as
\begin{align} \label{S2A}
\mathcal{S}_2 =  \sigma^2\left[\frac{k_{3,4}+1}{\Omega(k_{3,4}+1) - \Omega(k_{3,4})} + \frac{1}{\Omega(2)}\right],
\end{align}
Because $k_{3,4}$ are symmetric about $k=1$ (from the symmetry of $F$), the value of $\mathcal{S}_2$ is independent of the choice of $k_{3,4}$. Using the definition of the dispersion relation $\Omega$ \eqref{12}, \eqref{S2A} simplifies to
\begin{align}
\mathcal{S}_2 = \frac{\sigma^2}{2(1-5\beta)},
\end{align}
which is nonzero for $1/5 < \beta < 3/10$.  \begin{flushright} $\square$ \end{flushright}




\begin{thebibliography}{99} \footnotesize

 \bibitem{akers12}
B. Akers and W. Gao.
Wilton ripples in weakly nonlinear model equations.
\emph{Commun. Math. Sci.}
10(3): 1015-1024, 2012.

\bibitem{akersnicholls12}
B. Akers and D. P. Nicholls.
Spectral stability of deep two-dimensional gravity water waves: repeated eigenvalues.
\emph{SIAM J. Appl. Math.},
130(2): 81-107, 2012.

\bibitem{akersnicholls14}
B. Akers and D. P. Nicholls.
The spectrum of finite depth water waves.
\emph{European Journal of Mechanics-B/Fluids},
46: 181-189, 2014.

\bibitem{akers15}
B. Akers.
Modulational instabilities of periodic traveling waves in deep water.
\emph{Physica D},
300: 26-33, 2015

\bibitem{benjamin67}
T. B. Benjamin. Instability of periodic wave trains in nonlinear dispersive systems.
\emph{Proceedings, Royal Society of London, A,}
299:59-79, 1967.

\bibitem{benjaminfeir67}
T. B. Benjamin and J. E.  Feir. The disintegration of wave trains on deep water. part i. theory.
\emph{Journal of Fluid Mechanics},
27:417-430, 1967.

\bibitem{bridgesmielke95}
T. H. Bridges and A. Mielke. A proof of the Benjamin-Feir instability.
\emph{Archive for Rational Mechanics and Analysis},
133:145-198, 1995.

\bibitem{deconinckkapitula10}
B. Deconinck and T. Kapitula. The orbital stability of the cnoidal waves of the Korteweg-deVries equation. 
\emph{Phys. Letters A},
374: 4018-4022, 2010.

\bibitem{deconinck06}
B. Deconinck and J. N. Kutz. Computing spectra of linear operators using the Floque-Fourier-Hill method.
\emph{Journal of Computational Physics},
219(1): 296-321, 2006.

\bibitem{deconinck11}
B. Deconinck and K. Oliveras. The instability of periodic surface gravity waves.
\emph{Journal of Fluid Mechanics},
675:141-167, 2011.

\bibitem{deconinck17}
B. Deconinck and O. Trichtchenko.
High-frequency instabilities of small-amplitude of Hamiltonian PDE's.
\emph{Discrete \& Continuous Dynamical Systems-A},
37(3):1323-1358, 2017.

\bibitem{dias99}
F. Dias and C. Kharif.
Nonlinear gravity and gravity-capillary waves.
\emph{Annu. Rev. Fluid Mech.},
31: 331-341, 1999.

\bibitem{hammack93}
J. Hammack and D. Henderson.
Resonant interations among surface water waves.
\emph{Annu. Rev. Fluid Mech.},
25:55-97, 1993.

\bibitem{haragus06}
M. Haragus, E. Lombardi, and A. Scheel.
Spectral stability of wave trains in the Kawahara equation.
  \emph{Journal of Mathematical Fluid Mechanics},
 8(4):482-509, 2006.

\bibitem{haraguskapitula08}
M. Haragus and T. Kapitula.
On the spectra of periodic waves for infinite-dimensional Hamiltonian systems,
\emph{Phys. D},
237: 2649-2671, 2008.

\bibitem{haupt88}
S. E. Haupt and J. P. Boyd.
Modeling nonlinear resonance: a modification to the Stokes' perturbation expansion.
  \emph{Wave Motion},
 10(1):83-98, 1988.

\bibitem{johnson10}
M. A. Johnson, K. Zumbrun, and J. C. Bronski].]
On the modulation equations and stability of periodic
generalized Korteweg--de Vries waves via Bloch decompositions,
\emph{Phys. D},
239: 2067-2065, 2010.

\bibitem{kapitulapromislow13}
T. Kapitula and K. Promislow.
\emph{Spectral and Dynamical Stability of Nonlinear Waves.}
Springer, New York, 2013.

\bibitem{kato66}
T. Kato.
\emph{Perturbation Theory for Linear Operators.}
Springer-Verlag, Berlin, 1966.

\bibitem{kawahara72}
T. Kawahara. Oscillatory solitary waves in dispersive media.
\emph{J. Phys. Soc. Jpn.},
33: 1015-1024, 1972.

\bibitem{kollar19}
R. Koll\'{a}r, B. Deconinck, and O. Trichtchenko.
Direct characterization of spectral stability of small-amplitude periodic waves in scalar Hamiltonian problems via dispersion relation.
\emph{SIAM Journal on Mathematical Analysis},
51(4): 3145-3169, 2019.

\bibitem{krein51}
M. G. Krein. On the application of an algebraic proposition in the theory of matrices of monodromy.
\emph{Uspehi Matem. Nauk (N.S.)},
6(1(41)):171-177, 1951.


\bibitem{mackay86}
R. S. MacKay and P. G. Saffman.
Stability of water waves.
\emph{Proceedings of the Royal Society A},
406(1830):115-125, 1986.

\bibitem{nivala10}
M. Nivala and B. Deconinck. Periodic finite-genus solutions of the KdV equation are orbitally stable. 
\emph{Physica D}, 
239(13): 1147-1158, 2010.

\bibitem{pava14}
J. Pava and F. Natali. (Non)linear instability of periodic traveling waves: Klein-Gordon and KdV type equations.
\emph{Adv. Nonlinear Anal.},
3(2): 95-123, 2014.


\bibitem{stokes1847}
G. G. Stokes.
On the theory of oscillatory waves.
\emph{Trans. Camb. Phil. Soc.},
8:441-455, 1847.

\bibitem{trichtchenko16}
O. Trichtchenko, B. Deconinck, and J. Wilkening.
The instability of Wilton ripples.
\emph{Wave Motion},
66: 147-155, 2016.

\bibitem{trichtchenko18}
O. Trichtchenko, B. Deconinck, and R. Koll\'{a}r.
Stability of periodic traveling wave solutions to the Kawahara equation.
\emph{SIAM Journal on Applied Dynamical Systems},
17(4): 2761-2783, 2018.

\bibitem{whitham67}
G. B. Whitham. Non-linear dispersion of water waves.
\emph{Journal of Fluid Mechanics},
27:399-412, 1967.

\bibitem{wilton1915}
J. R. Wilton.
On ripples.
\emph{Phil. Mag. Ser.}
629(173): 688-700, 1915.




\end{thebibliography}
\end{document}